\documentstyle[11pt,bezier,figure]{article}

\begin{document}

\font\bbbld=msbm10 scaled\magstep1
\newcommand{\bfR}{\hbox{\bbbld R}}
\newcommand{\bfC}{\hbox{\bbbld C}}
\newcommand{\bfZ}{\hbox{\bbbld Z}}
\newcommand{\bfH}{\hbox{\bbbld H}}
\newcommand{\bfQ}{\hbox{\bbbld Q}}
\newcommand{\bfN}{\hbox{\bbbld N}}
\newcommand{\bfP}{\hbox{\bbbld P}}
\newcommand{\bfT}{\hbox{\bbbld T}}
\def\Sym{\mathop{\rm Sym}}
\newcommand{\halo}[1]{\Int(#1)}
\def\Int{\mathop{\rm Int}}
\def\Re{\mathop{\rm Re}}
\def\Im{\mathop{\rm Im}}
\newcommand{\union}{\cup}
\newcommand{\goesto}{\rightarrow}
\newcommand{\bdy}{\partial}
\newcommand{\n}{\noindent}
\newcommand{\p}{\hspace*{\parindent}}

\newtheorem{theorem}{Theorem}[section]
\newtheorem{assertion}{Assertion}[section]
\newtheorem{proposition}{Proposition}[section]
\newtheorem{lemma}{Lemma}[section]
\newtheorem{definition}{Definition}[section]
\newtheorem{claim}{Claim}[section]
\newtheorem{corollary}{Corollary}[section]
\newtheorem{observation}{Observation}[section]
\newtheorem{conjecture}{Conjecture}[section]
\newtheorem{question}{Question}[section]
\newtheorem{example}{Example}[section]

\newbox\qedbox
\setbox\qedbox=\hbox{$\Box$}
\newenvironment{proof}{\smallskip\noindent{\bf Proof.}\hskip \labelsep}%
                        {\hfill\penalty10000\copy\qedbox\par\medskip}
\newenvironment{remark}{\smallskip\noindent{\bf Remark.}\hskip \labelsep}%
                        {\hfill\penalty10000\copy\qedbox\par\medskip}
\newenvironment{remark1}{\smallskip\noindent{\bf Remark 1.}\hskip \labelsep}%
                        {\hfill\penalty10000\copy\qedbox\par\medskip}
\newenvironment{remark2}{\smallskip\noindent{\bf Remark 2.}\hskip \labelsep}%
                        {\hfill\penalty10000\copy\qedbox\par\medskip}
\newenvironment{proofspec}[1]%
                      {\smallskip\noindent{\bf Proof of Theorem 1.1.}
                        \hskip \labelsep}%
                        {\nobreak\hfill\hfill\nobreak\copy\qedbox\par\medskip}
\newenvironment{proofspec2}[1]%
                      {\smallskip\noindent{\bf Proof of Theorem 1.2.}
                        \hskip \labelsep}%
                        {\nobreak\hfill\hfill\nobreak\copy\qedbox\par\medskip}
\newenvironment{acknowledgements}{\smallskip\noindent{\bf Acknowledgements.}%
        \hskip\labelsep}{}

\setlength{\baselineskip}{1.0\baselineskip}

\title{Minimal Surfaces in $R^3$ with \\ 
Dihedral Symmetry.}
\author{Wayne Rossman} 
\date{}
\maketitle

\begin{abstract}
We construct new examples of immersed minimal surfaces with 
catenoid ends and finite total curvature, of both genus zero 
and higher genus.  In the genus zero case, we classify all such 
surfaces with at most $2n+1$ ends, and with symmetry group the natural
$\bfZ_2$ extension of the dihedral group $D_n$.
 \footnote{This research was 
supported by a fellowship from the Japan
Society for the Promotion of Science.}

The surfaces are constructed by proving 
existence of the conjugate surfaces.  We extend this method 
to cases where the conjugate surface of the fundamental piece 
is noncompact and is not a graph over a convex plane domain.  
\end{abstract}

\section{Introduction}

\begin{figure}
        \hspace{0.01in}
        \epsfxsize=1.4in
        \epsffile{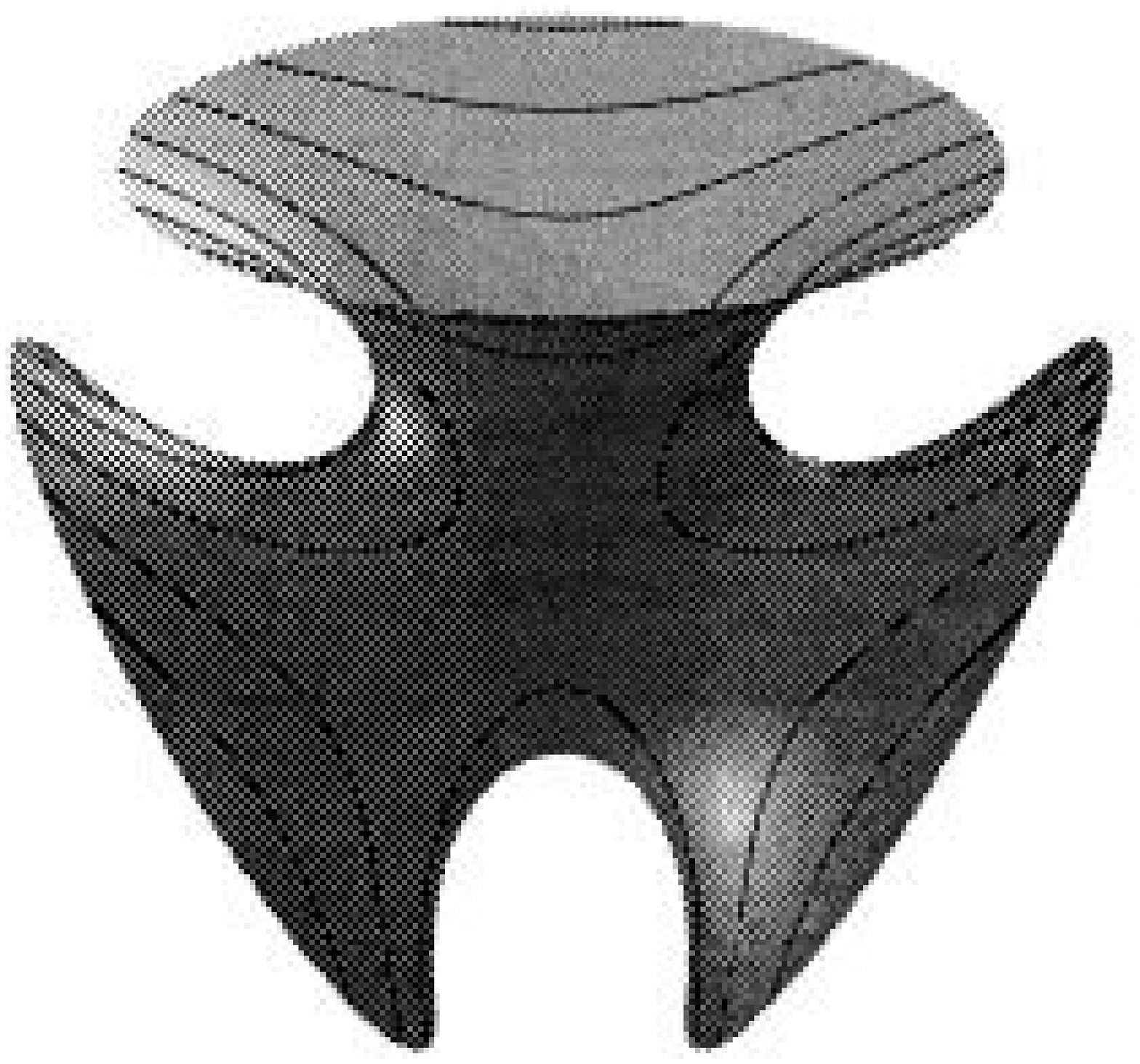}
        \hspace{0.01in}
        \epsfxsize=1.4in
        \epsffile{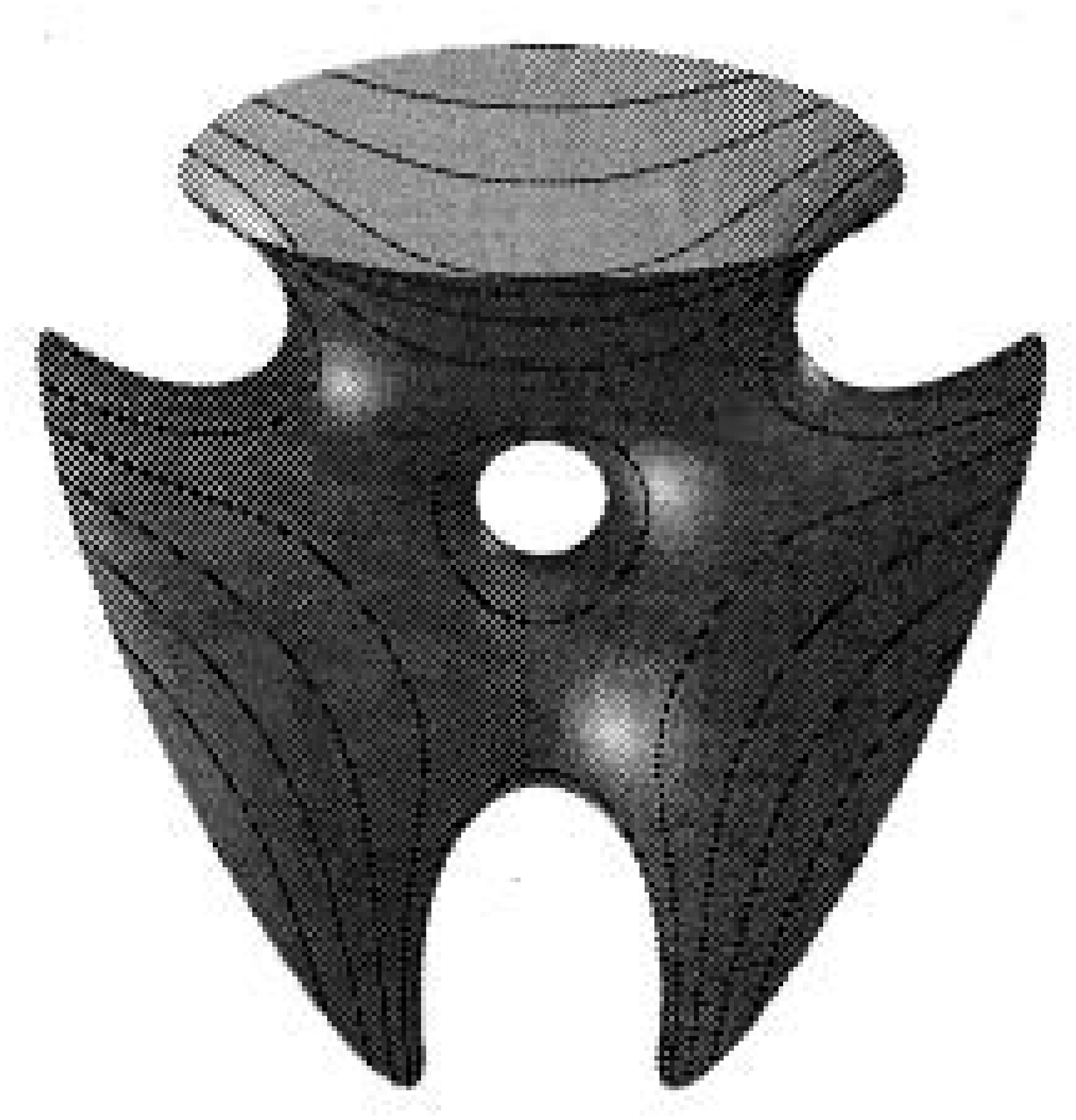}
        \hspace{0.01in}
        \epsfxsize=1.4in
        \epsffile{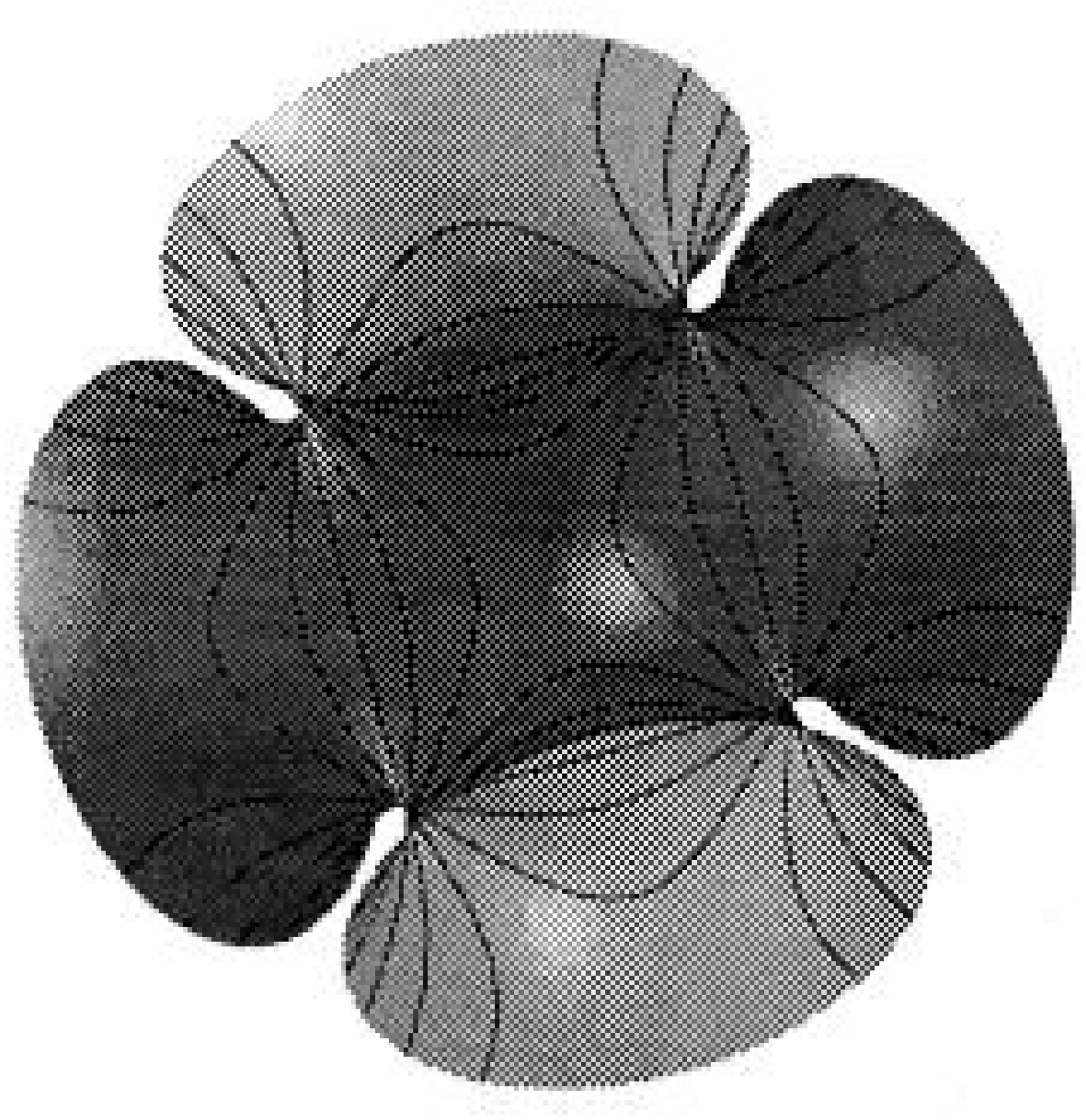}
        \hspace{0.01in}
        \epsfxsize=1.4in
        \epsffile{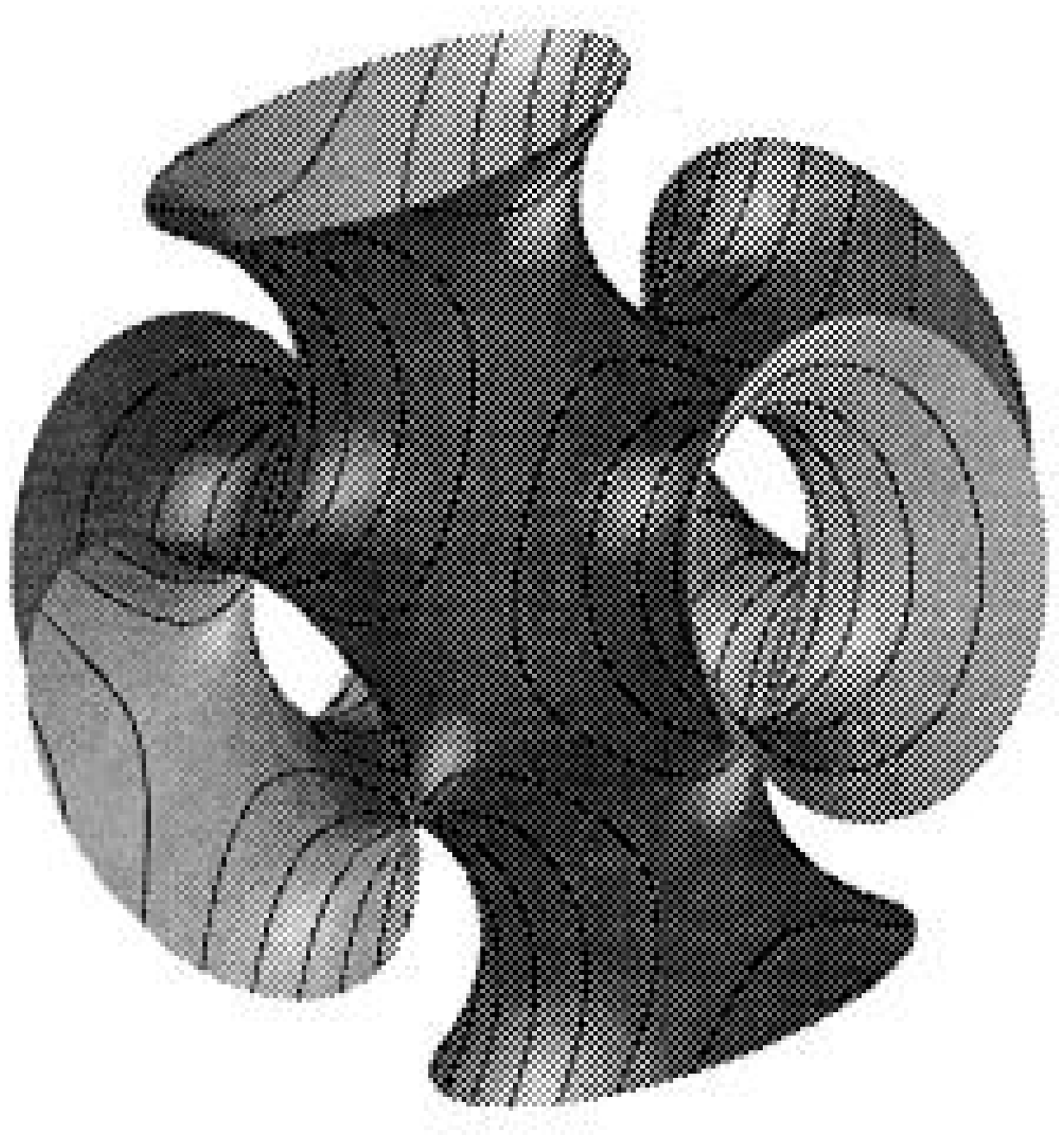}
        \hfill
        \caption{Jorge-Meeks 3-oid, genus-one 3-oid, genus-zero tetroid, genus-zero octoid}
        \label{old}
\end{figure}

\begin{figure}
        \hspace{0.01in}
        \epsfxsize=1.4in
        \epsffile{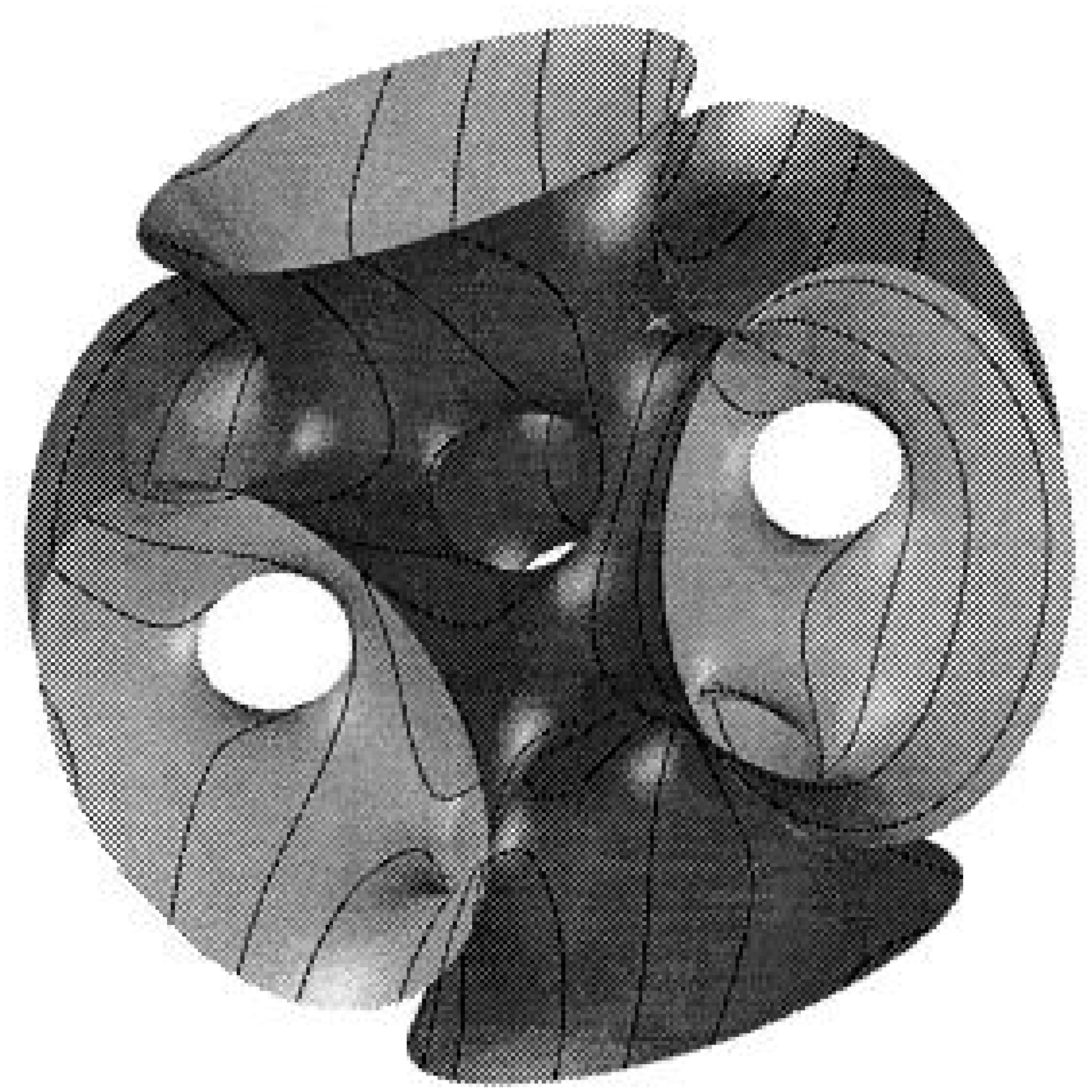}
        \hspace{0.01in}
        \epsfxsize=1.4in
        \epsffile{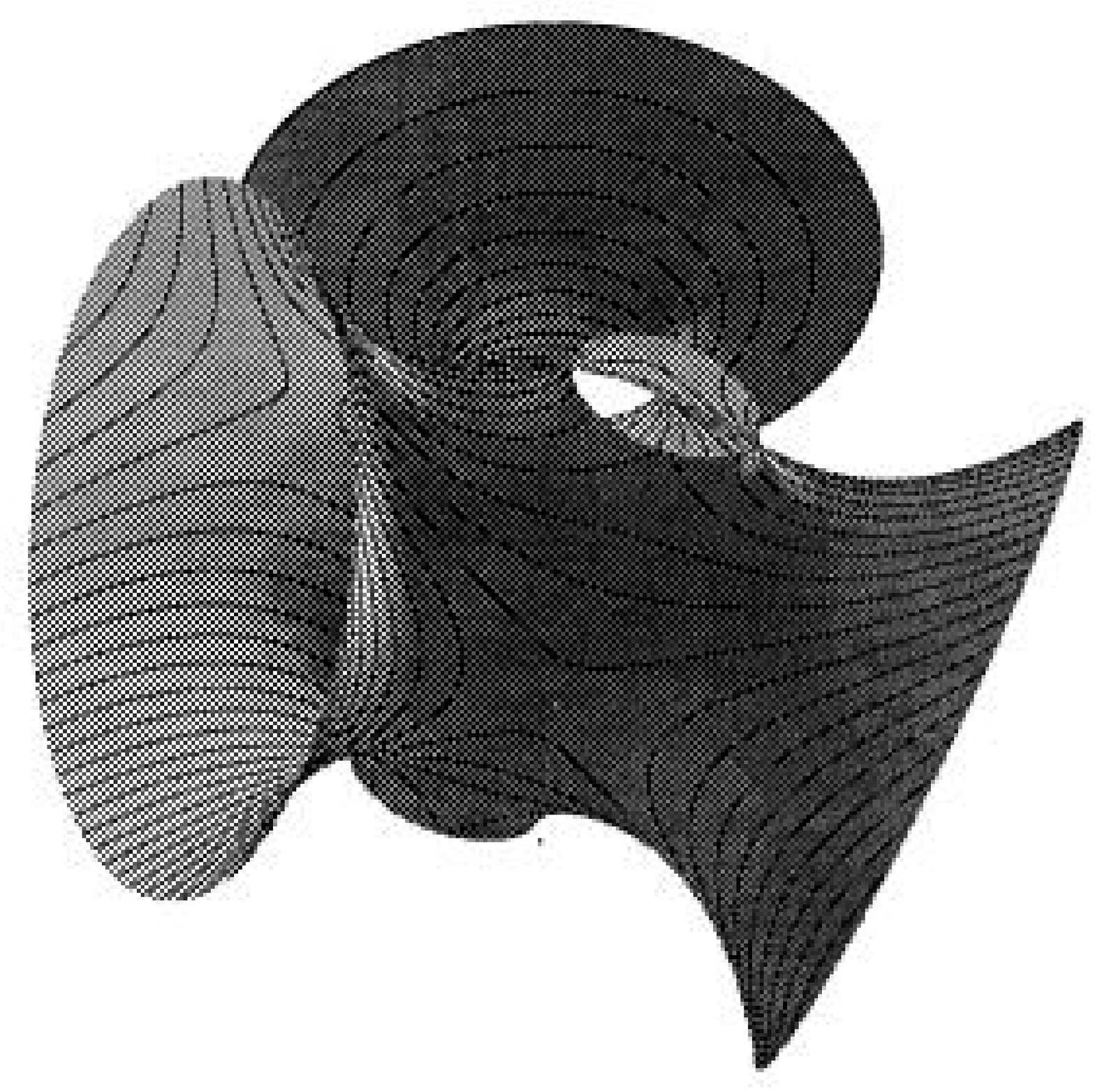}
        \hspace{0.01in}
        \epsfxsize=1.4in
        \epsffile{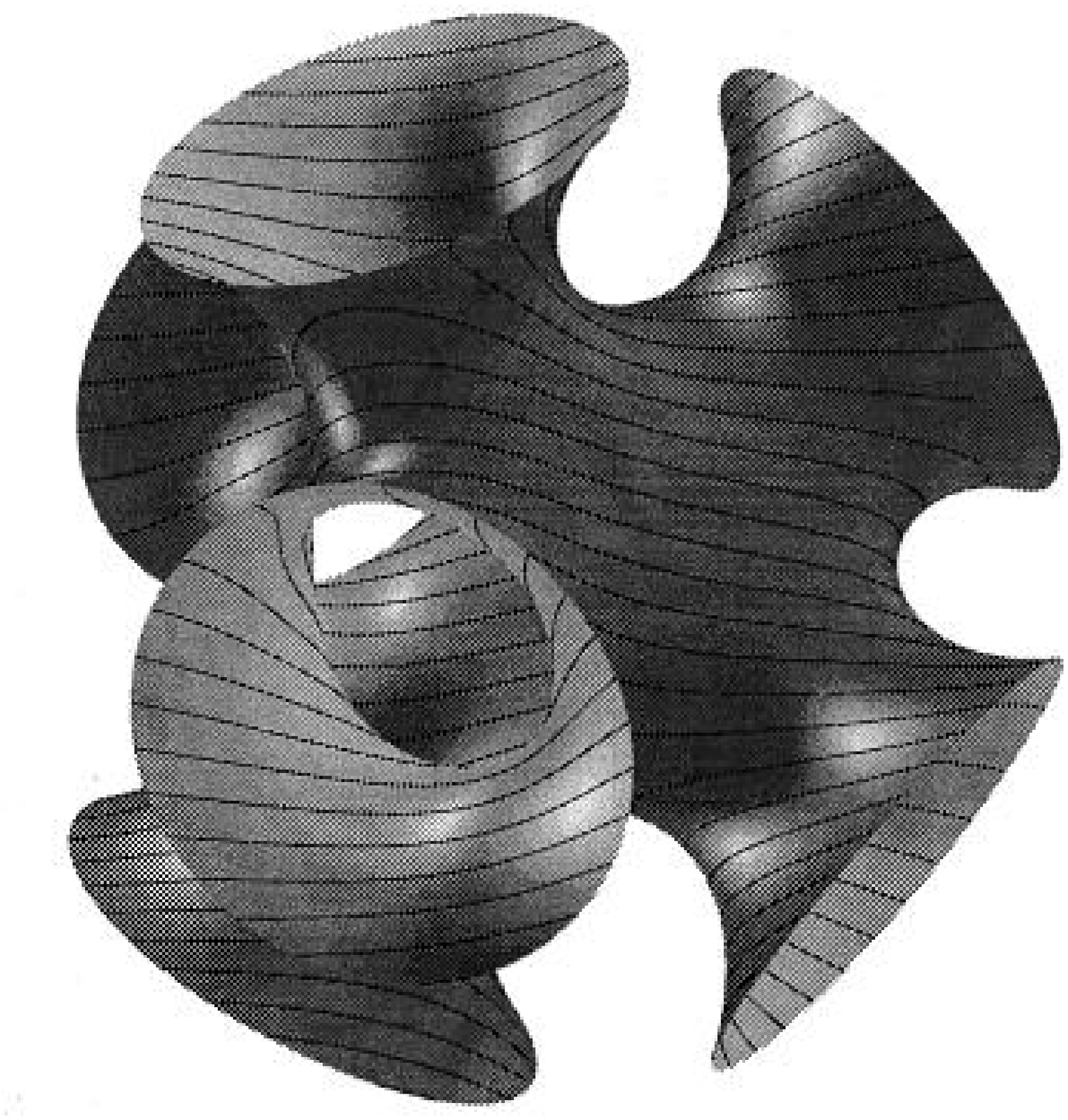}
        \hspace{0.01in}
        \epsfxsize=1.4in
        \epsffile{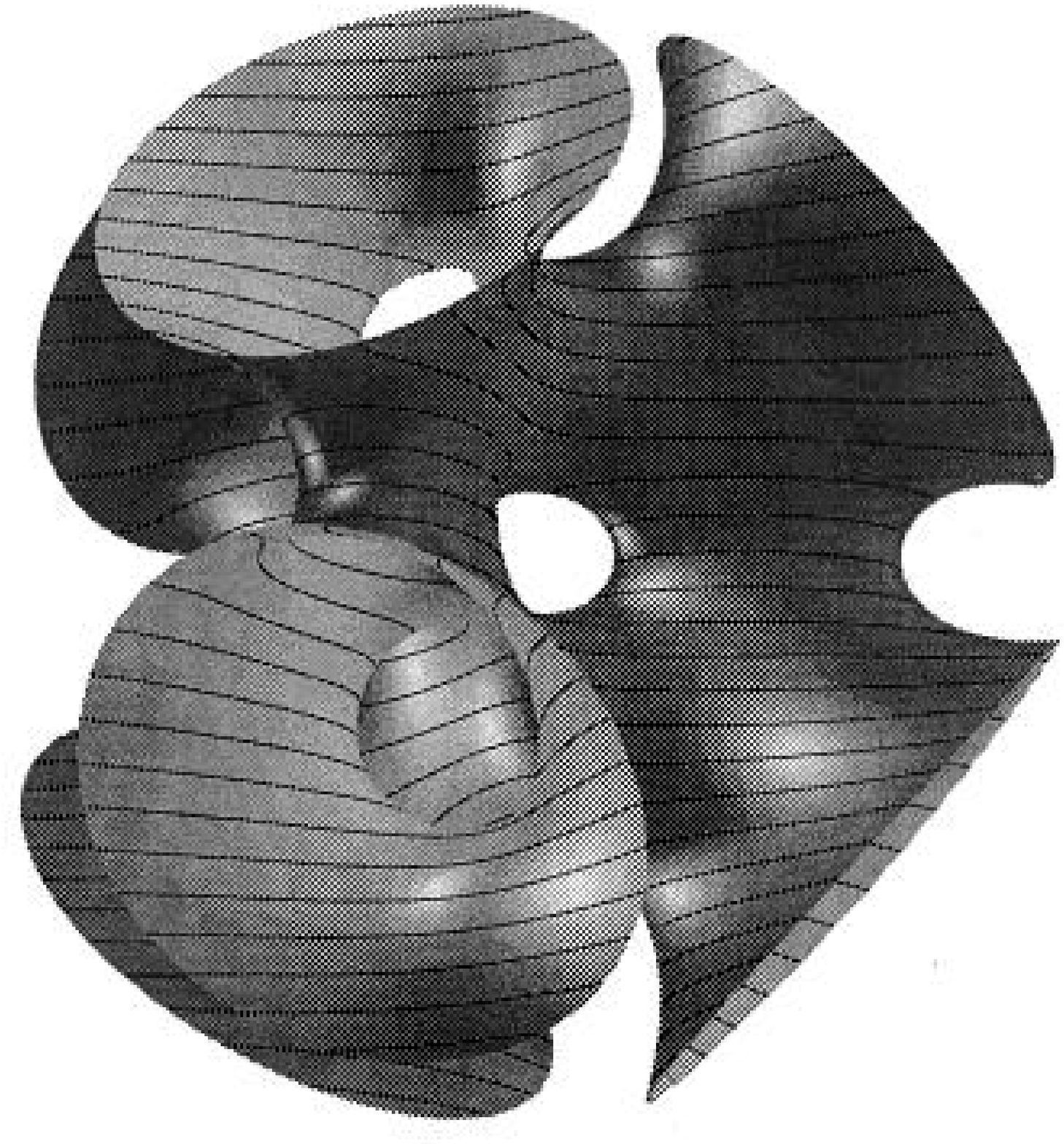}
        \hfill
        \caption{genus-seven octoid, Jorge-Meeks $n$-oid fence with $n = 3$, 
            ${\cal P}_0(2n,\theta)$ for $n=3$ and $\theta=45$ degrees, 
            ${\cal P}_{n-1}(2n,\theta)$ for $n=3$ and $\theta=45$ degrees}
        \label{between}
\end{figure}

Recently, new examples of immersed minimal surfaces of finite total 
curvature with catenoid ends have been found.  
Among these examples are: the genus-zero Jorge-Meeks
$n$-oid with symmetry group $D_n \times \bfZ_2$ \cite{JoMe}, the 
genus-one $n$-oid with symmetry group $D_n \times \bfZ_2$ \cite{BeRo}, the 
genus-zero Platonoids 
with symmetry groups isomorphic to the symmetry groups of the 
Platonic solids \cite{Xu}, \cite{Kat}, \cite{UmYa}, and 
the higher genus Platonoids 
with Platonic symmetry groups \cite{BeRo}.  (See Figures \ref{old} 
(1) - (4), \ref{between} (1).)  By $D_n \times \bfZ_2$, 
we mean the natural $\bfZ_2$ extension into $O(3)$ of the dihedral 
group $D_n \subset SO(3)$.  

\begin{figure}
        \hspace{0.2in}
        \epsfxsize=1.5in
        \epsffile{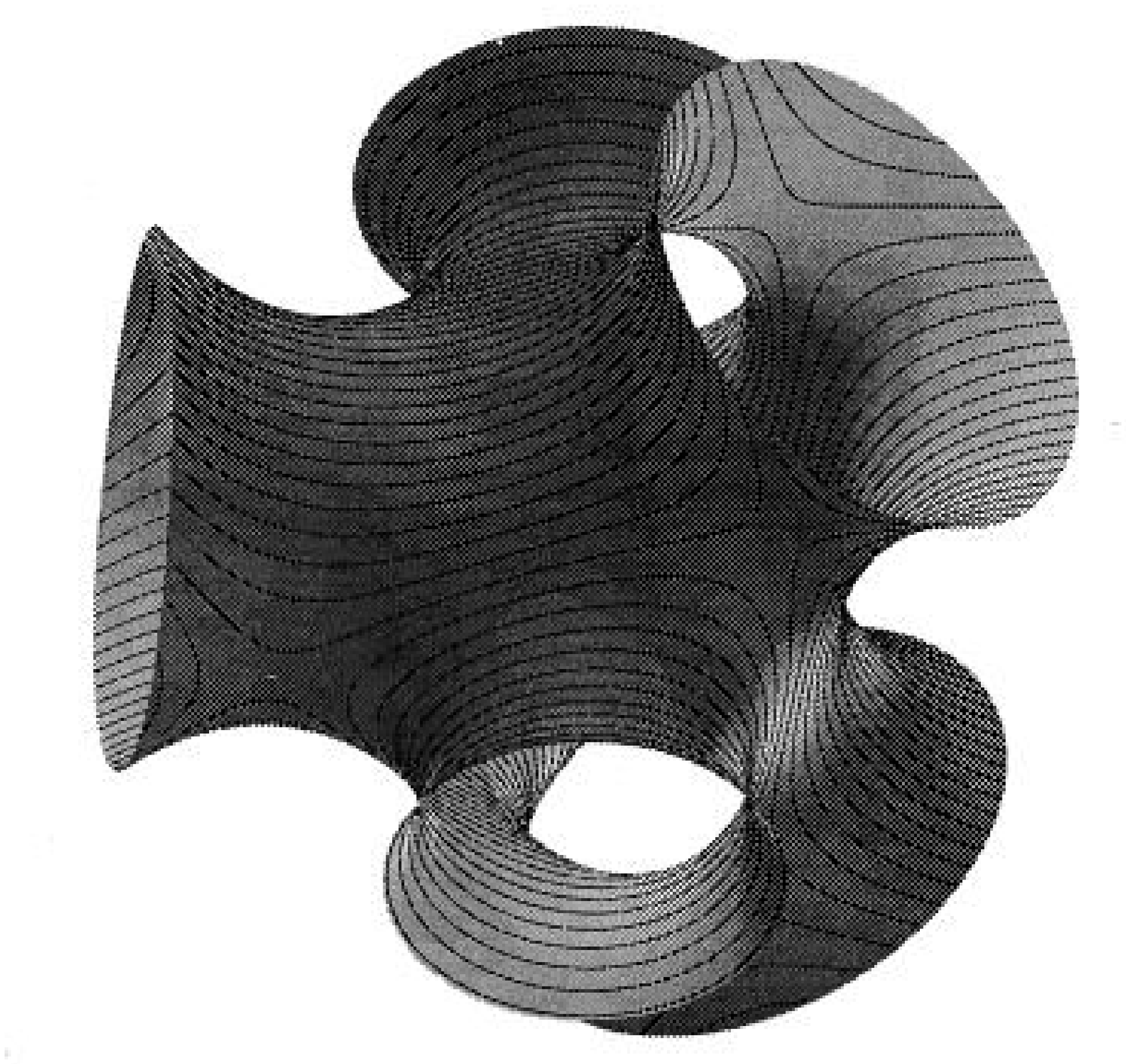}
        \hspace{0.2in}
        \epsfxsize=1.5in
        \epsffile{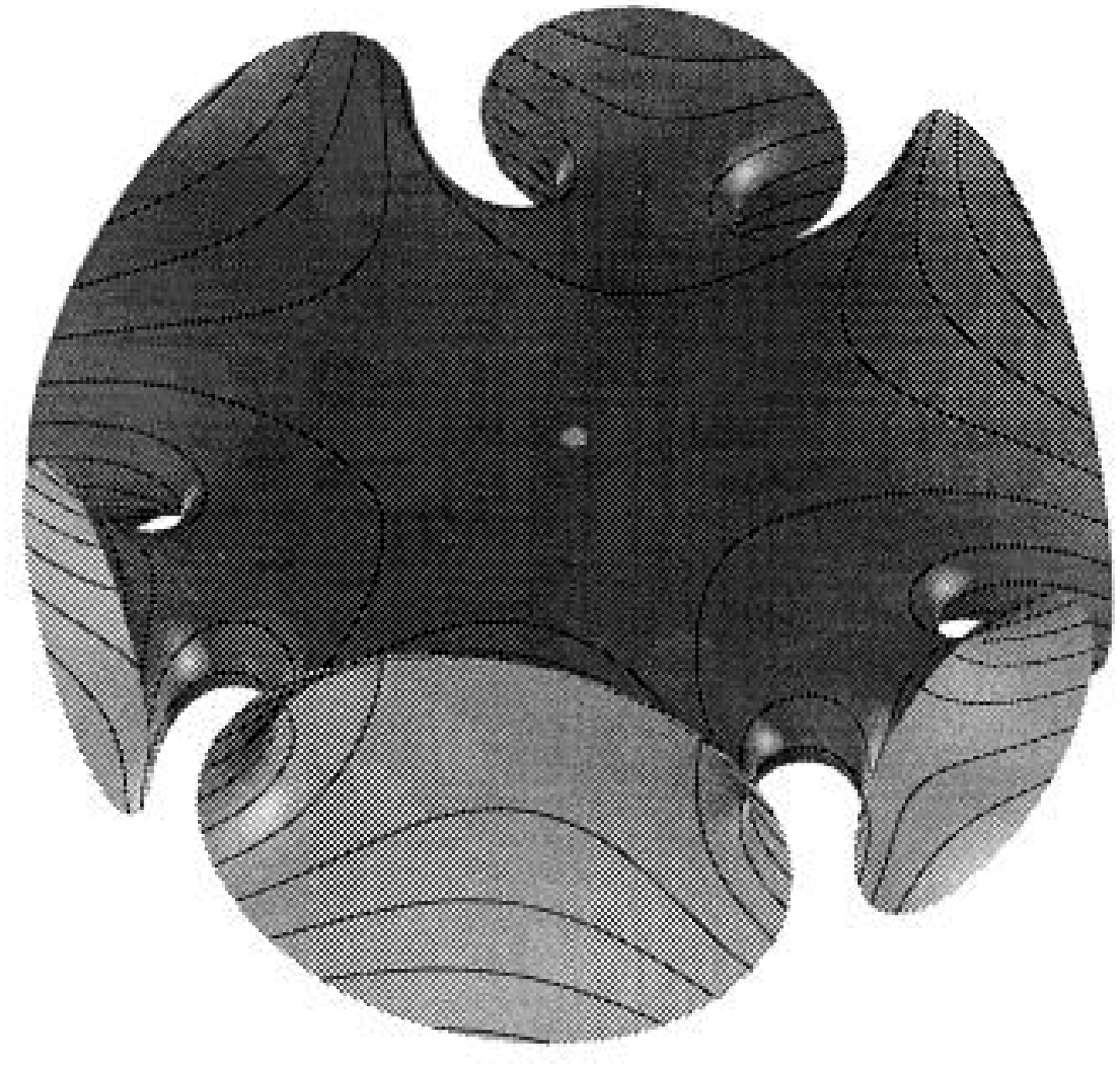}
        \hspace{0.2in}
        \epsfxsize=1.5in
        \epsffile{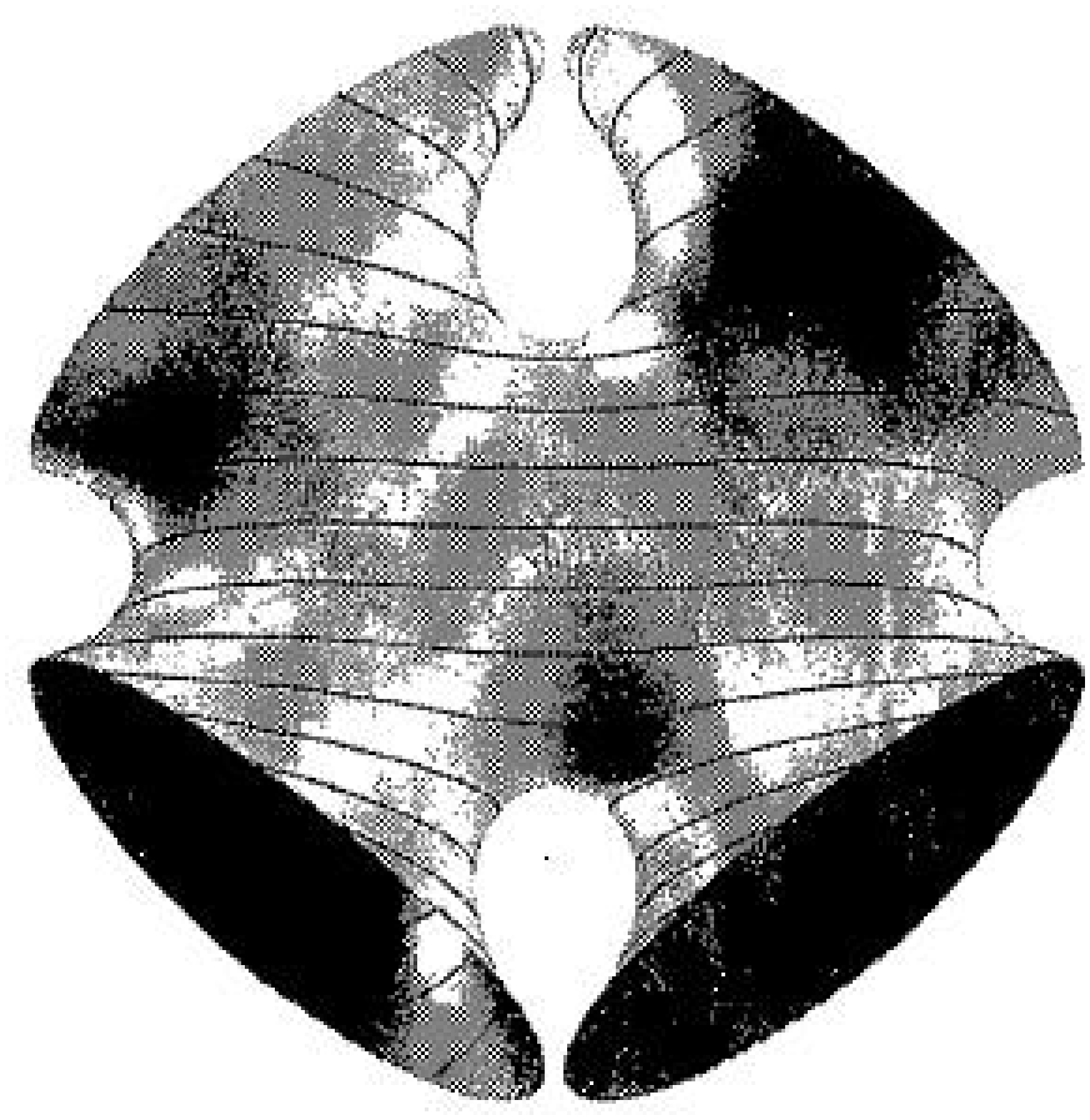}
        \hfill
        \caption{${\cal J}{\cal M}{\cal V}_0(n+2,w)$ for $n=3$ and $w \approx 5/2$, 
        ${\cal A}{\cal W}_0(2n,w)$ for $n=3$ and $w \approx 5/2$, 
        ${\cal A}{\cal A}_0(2n,\theta)$ for $n=2$ and $\theta=30$ degrees}
        \label{new}
\end{figure}

In this present work we find more examples with symmetry group $D_n \times 
\bfZ_2$ (see Figures \ref{between} (2) - (4), \ref{new} (1) - (3)), of 
both genus zero and higher genus.  Then, in the genus zero case, 
we classify all such surfaces that have at most $2n+1$ ends.  

To prove existence of these surfaces we use the conjugate surface 
construction, by an approach similar to that of 
\cite{BeRo}. Generally speaking, the conjugate surface construction 
seems to require a high degree of symmetry of the surface.
In fact, all of the known techniques for creating examples of minimal 
surfaces seem to benefit from symmetry assumptions.  

The examples we construct here are less symmetric than the 
examples mentioned in the first paragraph, 
in the sense that their fundamental 
pieces have higher total Gaussian curvature.  
It is therefore harder to prove existence of the 
conjugate surfaces to these fundamental pieces.  
Hence the constructions we need are more 
delicate than those in \cite{BeRo}.  For less 
symmetric surfaces, the conjugates may no longer lie over convex 
plane domains, thus making Nitsche's theorem no longer 
applicable.   In fact, they may not even be 
graphs, or may not even be embedded.  

Another consideration is the period problem.  In general, 
integration of the Weierstrass integral (described in the next 
section) about a nontrivial cycle produces a period vector.  
Translation by this period vector in $\bfR^3$ will produce an isometry 
of the surface.  Thus, if the surface has finite total 
curvature, the period vectors must be the zero vector for all 
cycles.  Ensuring that this is the case is called 
``removing'' or ``killing'' the periods.  

Often it is possible to remove one period with an 
intermediate-value-theorem argument, 
but for surfaces with more than one period cycle this will not 
suffice.  If there are two periods to remove, it is usually quite difficult 
to theoretically argue that both can be removed simultaneously.  
In our case we are able to make an argument to kill two 
periods (the proof of Theorem 1.2).  

Since it is well known \cite{Os} that the Gauss map of a 
finite-total-curvature minimal surface extends continuously to each 
end of the surface, we shall simply 
refer to the extended Gauss map at an end as the normal vector at that end.

The naming scheme for the surfaces discussed here is 
$\diamondsuit_A(B,C)$, where: $A$ is the genus of the surface, $B$ is 
the number of ends of the surface, and $C$ represents the parameter
for a one-parameter family of surfaces.  $C$ is either an angle
$\theta$ or a weight $w$ of an end, and is omitted when there is no
relevant one-parameter family.  In each case $\diamondsuit$ is
replaced by hopefully informative lettering.  (${\cal P}$ represents 
``prismoid'', ${\cal J}{\cal M}$ represents ``Jorge-Meeks surface'',
${\cal J}{\cal M}{\cal V}$ represents ``Jorge-Meeks surface with 
added vertical ends'', ${\cal A}{\cal W}$ represents ``surfaces with 
alternately weighted ends'', and ${\cal A}{\cal A}$ represents
``surfaces with alternating angles between ends''.) 

\begin{theorem}
{\sf (The prismoids)}  
For each $n \geq 2$, there exists a one-parameter family of 
immersed minimal surfaces ${\cal P}_0(2n,\theta), \, 0 < 
\theta < \pi/2$, satisfying the following:
\newcounter{num}
\begin{list}%
{\arabic{num})}{\usecounter{num}\setlength{\rightmargin}{\leftmargin}}
\item ${\cal P}_0(2n,\theta)$ has genus zero.
\item ${\cal P}_0(2n,\theta)$ has $2n$ catenoid ends, and the normal 
vector at each end makes an angle $\theta$ with a horizontal plane.
\item The symmetry group of ${\cal P}_0(2n,\theta)$ is $D_n \times \bfZ_2$.  
\end{list}
\end{theorem}

\begin{theorem}
{\sf (The higher genus prismoids)}  
For each $n \geq 2$, there exists a one-parameter family of 
immersed minimal surfaces ${\cal P}_{n-1}(2n,\theta), \, 0 < 
\theta < \pi/2$, satisfying the following:
\begin{list}%
{\arabic{num})}{\usecounter{num}\setlength{\rightmargin}{\leftmargin}}
\item ${\cal P}_{n-1}(2n,\theta)$ has genus $n-1$.
\item ${\cal P}_{n-1}(2n,\theta)$ has $2n$ catenoid ends, and the normal 
vector at each end makes an angle $\theta$ with a horizontal plane.
\item The symmetry group of ${\cal P}_{n-1}(2n,\theta)$ 
is $D_n \times \bfZ_2$.  
\end{list}
\end{theorem}

\begin{theorem}
{\sf (The genus-zero $n$-oids plus two ends)}  
For each $n \geq 2$, there exists a positive constant $c(n)$ 
so that, for any $w \geq c(n)$, there exists an 
immersed minimal surface 
${\cal J}{\cal M}{\cal V}_0(n+2,w)$ satisfying the following conditions:  
\begin{list}%
{\arabic{num})}{\usecounter{num}\setlength{\rightmargin}{\leftmargin}}
\item ${\cal J}{\cal M}{\cal V}_0(n+2,w)$ has genus zero.
\item ${\cal J}{\cal M}{\cal V}_0(n+2,w)$ 
has $n$ catenoid ends with weight one, and 
the normal vectors at these ends all lie within a 
horizontal plane and are symmetrically placed.  
\item ${\cal J}{\cal M}{\cal V}_0(n+2,w)$ 
has two catenoid ends of weight $w$ with 
vertical normals pointing in opposite directions.  
\item The symmetry group of ${\cal J}{\cal M}{\cal V}_0(n+2,w)$ 
is $D_n \times \bfZ_2$.  
\end{list}
\end{theorem}

By ``symmetrically placed'' in the second condition above, 
we mean that, up to a rotation of 
${\cal J}{\cal M}{\cal V}_0(n+2,w)$ if necessary, 
the $n$ ends with weight one have normal vectors whose 
stereographic projections to the complex plane are the $n$-th 
roots of unity.  This arrangement of ends is the same as for 
the Jorge-Meeks surface.  The remaining two ends of 
${\cal J}{\cal M}{\cal V}_0(n+2,w)$ 
have normal vectors whose stereographic projections are 
$z = 0$ and $z = \infty$.  

Roughly speaking, the weight of a catenoid end is a measure of 
the size of the catenoid to which it is asymptotic.  We give 
an exact definition in the next section.  In the previous theorem, 
the condition $w \geq c(n)$ 
seems to be unnecessarily restrictive, but is 
necessary for the proof we give here 
\cite{Xu}, \cite{KUY}, \cite{Kat}.

\begin{theorem}
{\sf (The $2n$-oids with alternating weights at the ends)}  
For each $n \geq 2$ there exists a one-parameter family of 
immersed minimal surfaces ${\cal A}{\cal W}_0(2n,w), \, 0 < w < \infty$, 
satisfying the following: 
\begin{list}%
{\arabic{num})}{\usecounter{num}\setlength{\rightmargin}{\leftmargin}}
\item ${\cal A}{\cal W}_0(2n,w)$ has genus zero.
\item ${\cal A}{\cal W}_0(2n,w)$ has $2n$ catenoid ends.  These ends 
have normal vectors all lying within a common plane and 
symmetrically placed.  (That is, up to a rotation of the 
surface if necessary, these normal vectors stereographically project to the 
$2n$-th roots of unity.)  
\item Half of the $2n$ ends have weight one, and the other $n$ 
ends have weight $w$, and they alternate between each other.  
\item The symmetry group of 
${\cal A}{\cal W}_0(2n,w)$ is $D_n \times \bfZ_2$.
\end{list}
\end{theorem}

\begin{theorem}
{\sf (The $2n$-oids with alternating angles between the ends)}  
For each $n \geq 2$ there exists a one-parameter family of 
immersed minimal surfaces 
${\cal A}{\cal A}_0(2n,\theta), \, 0 < \theta < 
\pi/n$, satisfying the following: 
\begin{list}%
{\arabic{num})}{\usecounter{num}\setlength{\rightmargin}{\leftmargin}}
\item ${\cal A}{\cal A}_0(2n,\theta)$ has genus zero.
\item ${\cal A}{\cal A}_0(2n,\theta)$ 
has $2n$ catenoid ends.  These ends all 
have weight one, and 
have normal vectors all lying within a common plane.  
\item The angles between adjacent ends alternate between $\theta$ and 
$(2\pi - n\theta)/n$.  
\item The symmetry group of 
${\cal A}{\cal A}_0(2n,\theta)$ is $D_n \times \bfZ_2$.
\end{list}
\end{theorem}

\begin{theorem}
{\sf (Classification)}  
Any genus-zero catenoid-ended immersed minimal surface 
with symmetry group $D_n \times \bfZ_2$ and at most $2n+1$ ends 
is either the Jorge-Meeks $n$-oid, 
${\cal P}_0(2n,\theta)$ for some $\theta \in \bfR$, 
${\cal J}{\cal M}{\cal V}_0(n+2,w)$ for some
$w \in \bfR$, ${\cal A}{\cal W}_0(2n,w)$ for some $w \in \bfR$, or 
${\cal A}{\cal A}_0(2n,\theta)$ for some $\theta \in \bfR$.  
\end{theorem}

In this classification we do not need to place any restrictions on the
ranges of $\theta$ or $w$.  The proof of this theorem given 
in Section 4 is independent of the values of $\theta$ and $w$.  (Note
that ${\cal J}{\cal M}{\cal V}_0(n+2,w)$ 
with $w=0$ and ${\cal A}{\cal W}_0(2n,w)$ with $w=0$ are simply
the Jorge-Meeks $n$-oid.)  
We caution, however, that we have not proven existence of 
${\cal P}_0(2n,\theta)$ (resp. 
${\cal A}{\cal A}_0(2n,\theta)$) 
when $\theta \not\in (0,\pi/2)$ (resp. $\theta \not\in (0,\pi/n)$), 
and of ${\cal J}{\cal M}{\cal V}_0(n+2,w)$ 
(resp. ${\cal A}{\cal W}_0(2n,w)$) when $w<c(n)$ (resp. $w<0$).  

The author wishes to thank: Jorgen Berglund, Frank Morgan, Shin 
Nayatani, Kotaro Yamada for helpful conversations; Rob Kusner for 
suggesting the research and for critical readings; Shin Kato, Masaaki 
Umehara for helpful conversations and for correcting an error in a 
preliminary draft; and Seiki Nishikawa for his support.  

\section{Preliminaries}
\subsection{The conjugate surface construction}

Consider a simply-connected finite-total-curvature immersed minimal 
surface $M$ in $\bfR^3$ with a boundary consisting of a 
finite number of piece-wise smooth curves.  
As proven by Enneper and Weierstrass, there exists a meromorphic function 
$g$ and a holomorphic 1-form $\eta$ defined 
on the unit disk in the complex plane such that $M$ has the 
parametrization
\[ \Phi(p) = \mbox{Re}\int_{p_0}^{p}
\; \left( \begin{array}{c}
        (1-g^2)\eta \\
        i(1+g^2)\eta \\
        2g\eta
        \end{array}
\right) \; , \; \; \; p \in \{z \in \bfC \; \, \mbox{such that} 
\; \, |z| \leq 1\} \; \; \; . \]

We refer to $\{g,\eta\}$ as the Weierstrass data for $M$, and to
$\Phi$
as the Weierstrass representation of $M$.  The map 
$g$ is stereographic projection of the Gauss map from the sphere 
to the complex plane.  
The conjugate surface $M_{\mbox{\footnotesize conj}}$ 
of $M$ is the minimal surface with 
the same parametrization, but with Weierstrass data $\{g,i\eta\}$; 
that is, $\eta$ is replaced with $i\eta$ in the parametrization 
above.  We shall call this conjugate Weierstrass representation 
$\Phi_{\mbox{\footnotesize conj}}(p)$.  Note that the conjugate of 
the conjugate of $M$ is given by the Weierstrass data $\{g,-\eta\}$,
giving us the original surface reflected through the origin.

Thus we have the maps $z \longmapsto \Phi(z)$ and $z \longmapsto
\Phi_{\mbox{\footnotesize conj}}(z)$ from the unit disk to $M$ and 
$M_{\mbox{\footnotesize conj}}$, respectively.  This induces a map 
$\phi = \Phi_{\mbox{\footnotesize conj}} \circ \Phi^{-1}$, the {\em 
conjugate map}, from $M$ to $M_{\mbox{\footnotesize conj}}$.  
The conjugate map $\phi$
is an isometry and preserves the Gauss map.  It also has the 
following property, which we shall use in an essential way: 
$\phi$ maps planar principal curves in $M$ to 
planar asymptotic curves in $M_{\mbox{\footnotesize conj}}$, 
and maps planar asymptotic 
curves in $M$ to planar principal curves in $M_{\mbox{\footnotesize conj}}$. 
From this we can conclude that {\em $\phi$ maps 
non-straight planar geodesics to straight lines, and vice versa}.  
And since the Gauss map is preserved by $\phi$, it follows that if 
$\phi$ maps a non-straight 
planar geodesic $\alpha \subset M$ to a line segment 
$\beta \subset M_{\mbox{\footnotesize conj}}$, 
the line segment $\beta$ must be 
perpendicular to the plane containing $\alpha$.  

In the cases we consider here, $M$ is bounded by 
piecewise smooth boundary curves that consist of a finite 
number of planar geodesics, and hence $M_{\mbox{\footnotesize conj}}$ 
is bounded by 
piecewise smooth boundary curves that consist of a finite 
number of line segments, rays, and complete lines.  

Recall from \cite{Scn1} that an end of a complete minimal immersion 
in $\bfR^3$ is a {\em regular} end if a neighborhood of this end 
is a graph $f$ with bounded slope 
over some plane (without loss of generality, the 
$x_1x_2$-plane), so that $f$ has the following asymptotic 
behavior:  
\[ f(x_1,x_2) = a\log(\sqrt{x_1^2+x_2^2}) + b + 
\frac{c_1x_1 + c_2x_2}{x_1^2+x_2^2} + {\cal O}(
\frac{1}{x_1^2+x_2^2})   \; \; \; . \]
If $a=0$ we have a planar end, and if $a \neq 0$ we have a 
catenoid end.  In this paper we shall use the terminology 
more loosely.  We shall say that a minimal end is a catenoid 
end (resp. planar end) if it satisfies the above 
asymptotic condition with $a \neq 0$ (resp. $a=0$), even if the 
minimal immersion has a nonempty boundary and there exist 
boundary curves which extend to the end.  With this more 
general definition in mind, we shall say that a 
minimal surface has a {\em helicoid end} if
the corresponding end of the conjugate surface is a catenoid end.  

For more detailed information on the conjugate surface 
construction, see \cite{Ka1}, \cite{Ka2}, \cite{Ka3}, 
\cite{Ka4}, \cite{BeRo}.   

\subsection{Weights}

Here we define a useful quantity (see \cite{KKS}) that is a 
vector associated to each Jordan curve in a minimal surface.  
We then describe some properties of this weight vector 
that are pertinent to our situation.  

\begin{definition}
Let $\alpha_1,...,\alpha_k$ be the 
boundary curves of a compact immersed minimal 
surface $M \subset \bfR^3$.  
Let $\upsilon$ be the outward pointing unit conormal of $M$ along
$\alpha_i$.  Then the {\em weight} (also called {\em flux})
of the boundary curve $\alpha_i$ is
\[ w(\alpha_i) = \int_{\alpha_i}\upsilon \, ds \; . \]
\end{definition}

It is a well-known application of Stoke's theorem 
that the weight vectors satisfy a ``balancing'' condition 
\[ \sum_{i=1}^k w(\alpha_i) = 0 \; . \]  
Furthermore, the weight vector can be defined by the same 
integral for any closed curve $\alpha$ on any complete oriented 
minimal surface $M$, up to a sign.  The signature depends on the 
choice of orientation of the conormal along $\alpha$.  

It readily follows that 
the weight vector $w(\alpha)$ associated to each closed curve 
$\alpha \subset M$ is an invariant of homology.  That is to say, 
if $\alpha, \beta \subset M$ are homologous closed curves with 
the same orientation of the conormal, then $w(\alpha) = w(\beta)$ 
(cf. \cite{HoMe}, \cite{KKS}).  

Thus, by considering any closed loop about each end of a complete 
finite-total-curvature immersed 
minimal surface, there is a well defined weight vector 
associated to each end.  It is easily seen that an embedded 
finite-total-curvature end is 
a catenoid end if and only if it has a nonzero weight vector, and 
is a planar end if and only if its weight vector is zero.  

We say that a set of 
$n$ vectors is in a balanced configuration if their sum is 
zero.  If a minimal surface $M$ with $n$ catenoid ends has these $n$ 
vectors as the weight vectors of its ends, we say that this 
configuration of vectors is {\em realized} by $M$.  
Clearly, if the configuration is not balanced, it cannot be 
realized by any minimal surface.  

Suppose a catenoid end $E$ of a minimal surface $M$ is asymptotic to a 
catenoid $\cal C$, and $E$ has weight vector $w(E)$.  From 
the Weierstrass representation, we can see that $\cal C$ 
is, up to a rigid motion of $\bfR^3$, a catenoid whose Weierstrass data is 
\[ g = \frac{1}{z} \; \; , \; \; \; \eta = \frac{|w(E)|}{4\pi}dz \; \; , \] 
where the base Riemann manifold is $\bfC \setminus \{0\}$.  
It follows that $|w(E)|$ is proportional to the ``size'' of $E$.

We can then see that $|w(E)|$ is the length of the fundamental 
period vector of the helicoid which is conjugate to $\cal C$.  
(If one also considers periods of the helicoid that do not 
preserve orientation, then the length of the fundamental
period vector is $|w(E)|/2$.)  

\subsection{Known results}

Before proving the theorems in this paper, we give some preliminary 
results that will be needed for the proofs.  Here we 
state some results that come from previous works, and in 
Subsection 2.4 we prove three lemmas that are designed specifically 
for our purposes.  

The following well-known lemma is the maximum principle for minimal 
surfaces.  It is a special case of a lemma in \cite{Scn1}, 
and is proven there.  We apply this lemma later in a variety of 
situations.

\begin{lemma}
        {\em (Interior Maximum Principle)} Let $M_1$ and
$M_2$ be minimal
surfaces in $\bfR^3$.  Suppose $p$ is an interior
point of both $M_1$, $M_2$, and suppose $T_p(M_1) = T_p(M_2)$.
If $M_1$ lies to one side of $M_2$ near $p$, then $M_1 = M_2$.

        {\em (Boundary Point Maximum Principle)} Suppose $M_1$,
        $M_2$ have
$C^2$-boundaries $C_1$, $C_2$, respectively, and
suppose $p$ is a point of both
$C_1$, $C_2$.  Furthermore, suppose the tangent planes of both
$M_1$, $M_2$ and $C_1$, $C_2$ agree at $p$: that is to say, suppose
$T_p(M_1) = T_p(M_2)$, $T_p(C_1) =
T_p(C_2)$.  If, near $p$, $M_1$ lies to one side
of $M_2$, then $M_1 = M_2$.
\end{lemma}

The following theorems are special cases of a result by 
Meeks-Yau \cite{MeYa}, and a result by Nitsche \cite{Ni}, \cite{JeSe}, 
\cite{BeRo}, \cite{MeYa}.  These theorems will be needed later 
to show that the Plateau solutions for certain polygonal contours 
are embedded.  

\begin{theorem}
Let $\hat{M}$ be a $3$-dimensional compact submanifold of $\bfR^3$ so that 
$\partial \hat{M}$ is piecewise smooth and consists of the smooth pieces 
$\{H_1,...,H_l\}$.  Assume the following two conditions:
\begin{list}%
{\arabic{num})}{\usecounter{num}\setlength{\rightmargin}{\leftmargin}}
\item Each $H_i$ is a compact subset of some minimal surface in 
$\bfR^3$.
\item Whenever $H_i$ and $H_j$ meet along a curve, the angle between the 
two surfaces is at most 180 degrees, with respect to the region $\hat{M}$.  
\end{list}
Let $\alpha$ be a Jordan curve in $\partial \hat{M}$.  Then there exists a 
branched minimal immersion from a disk $D$ into $\hat{M}$ with boundary 
$\alpha$, which is smooth in the interior of $D$ and has minimal 
area among all such maps.  Furthermore, any branched minimal 
immersion of the above type must be an embedding.  
\end{theorem}

\begin{theorem}
Let $D$ be a bounded convex domain in a horizontal plane, so that its
boundary $\partial D$ is piecewise smooth.  Let $\partial
\tilde{D} =
\partial D \setminus \{p_1,...,p_r\}$.  Then there exists a solution 
(as a graph over $D$) of the 
minimal surface equation in $D$ taking on preassigned bounded
continuous
data on the arcs of $\partial \tilde{D}$.  As a surface, this solution
contains vertical
line segments over the jump discontinuities of the boundary data.
\end{theorem}

\subsection{Lemmas}

Consider a finite-topology minimal surface $M$ (with boundary 
$\partial M$) with an end that is a 180 degree 
arc of a helicoid end.  Denote a neighborhood of this end 
by $E$.  By rotating if necessary, we may 
assume that the normal vector at this end is vertical.  
Suppose that outside a compact ball in $\bfR^3$ the boundary $\partial 
E$ is a pair of straight (necessarily horizontal) 
rays $r_1, r_2$.  The conjugate surface $E_{\mbox{\footnotesize 
conj}}$ of $E$ is a 
surface with a 180 degree arc of a catenoid end 
that, outside a compact ball 
in $\bfR^3$, is bounded by two infinite planar geodesics 
$s_1, s_2$ asymptotic to catenaries.  The curves $s_1, s_2$ 
lie in parallel vertical planes.  For this situation, we 
have the following lemma.

\begin{lemma}
The two planar geodesics $s_1, s_2 \in \partial E_{\mbox{\footnotesize conj}}$ 
lie in the same plane if and only if the two conjugate straight boundary 
rays $r_1, r_2 \in \partial E$ lie in a common vertical plane.  
\end{lemma}

\begin{proof}
Assume that $r_1$ and $r_2$ lie in a common vertical plane.  
Let Rot$_1$ be the 180 degree rotation about the line containing 
$r_1$, and let Rot$_2$ be the 180 degree rotation about the line 
containing $r_2$.  

The surface $E \cup$Rot$_2(E)$ is a smooth embedded end 
asymptotic to a 360 degree arc of a helicoid end; and outside 
of a compact ball in $\bfR^3$, it is bounded by two parallel 
rays $r_1$, Rot$_2(r_1)$, which also lie in a common vertical plane.  

We choose the orientation for $E \cup$Rot$_2(E)$ so that the normal 
vector at the end is $(0,0,-1)$.  Thus 
the Weierstrass data for this end can be 
given, in a punctured neighborhood of the origin in $\bfC$, as 
\[ g = c_1 z + {\cal O}(z^2) \; \; , \; \; 
\eta = (\frac{c_2}{z^2} + \frac{c_3}{z} + 
{\cal O}(1)) \, dz \; \; \; .  \]
The conformal transformation $z \longmapsto z/c_1$ preserves 
the origin, so we may therefore assume that $c_1 = 1$.

Since Rot$_2 \circ$Rot$_1$ is a vertical translation, the surface 
$E \cup$Rot$_2(E)$ is a portion of a helicoid end that is 
periodic in the $x_3$ direction.  
Therefore, in the Weierstrass representation with 
this data, integrating around 
a small circle $\{z \in \bfC \; \mbox{such that} \, |z| = \epsilon\}$ 
about the origin results in a vertical period.  
From an examination of 
the first two coordinates of the Weierstrass representation, we 
see that $c_3$ must be 0.  

Now consider $E_{\mbox{\footnotesize conj}}$ and its reflection 
Ref($E_{\mbox{\footnotesize conj}}$) 
through the plane containing $s_2$.  We wish to conclude that 
Ref($s_1$) and $s_1$ are the same curve.  
The Weierstrass data for this catenoid end is 
\[ g = z + {\cal O}(z^2) \; \; , \; \; 
\eta = (\frac{ic_2}{z^2} + {\cal O}(1)) \, dz \; \; \; .  \]
Since this is a catenoid end with vertical normal vector, it 
satisfies the asymptotic condition in Subsection 2.1, therefore it 
cannot have any periodicity in the $x_3$ direction.  
Examining the third coordinate of the 
Weierstrass representation for this data shows 
that $c_2$ is purely imaginary.  It follows that integrating 
around $\{z \in \bfC \; \mbox{such that} \, |z| = \epsilon\}$
produces the zero vector.  
Thus Ref($s_1$) and $s_1$ are indeed the same curve.  
Hence $s_1$ and $s_2$ lie in the same plane.  

The above argument can be reversed to produce the converse 
conclusion.  
\end{proof}

The following lemma will be needed later to extend compact 
embedded Plateau solutions to stable noncompact embedded minimal 
surfaces.  We use the term stable in the following sense:
A noncompact minimal surface $M$ 
(possibly with boundary) is {\em stable} if the
second derivative of area is nonnegative at $M$ for all
smooth variations of the surface with compact support (and fixing 
the boundary $\partial M$).  

\begin{lemma}
Let $\{C_i\}_{i=1}^\infty$ be a sequence of compact Jordan contours 
in $\bfR^3$ so that the following conditions hold:
\begin{list}%
{\arabic{num})}{\usecounter{num}\setlength{\rightmargin}{\leftmargin}}
\item Each $C_i$ is a piecewise smooth contour consisting of a 
finite number of line segments.  
\item Each $C_i$ bounds a least-area minimal disk $M_i$.  
\item For any ball $B_R$ of radius $R$ in $\bfR^3$, 
there exists $N_R \in \bfZ$ so that $C_i \cap B_R = C_j \cap B_R$ 
for any $i, j \geq N_R$.  
\item There exists a fixed $\delta > 0$ and 
a compact $3$-dimensional region $\hat{M} \subset B_{1/\delta} 
\subset \bfR^3$ so that $M_i \cap \hat{M} \neq 
\phi$ and ${\rm dist}(C_i,\hat{M}) > \delta \,$ for all $i$.  
\item $\{C_i\}_{i=1}^\infty$ converges (in the topology of compact uniform 
convergence) to a noncompact contour $C$, and 
$C$ is a piecewise smooth contour consisting of a finite 
number of line segments, rays, and complete lines.
\end{list}
Then a subsequence of $\{M_i\}_{i=1}^\infty$ converges to a nonempty 
stable minimal surface $M$ (possibly disconnected) with boundary 
$C$.  Furthermore, if each $M_i$ is embedded, then $M$ is 
embedded.  
\end{lemma}

\begin{proof}
Schoen \cite{Scn2} has proven that the Gaussian curvature on 
a stable minimal surface $M \subset \bfR^3$ is bounded by 
$|K(p)| \leq c/{r^2}$, where $r$ is the distance within $M$ from the 
point $p \in M$ to the boundary $\partial M$, and 
$c$ is a universal constant.  
Let ${\cal N}_\epsilon(C)$ be an $\epsilon$-neighborhood of $C$.  
From Schoen's estimate and the fact that $C_i \cap B_R = C \cap B_R$ 
for $i$ large enough, we see that the function $|K|$ 
is bounded by $c/{\epsilon^2}$ on $M_i \cap (B_R \setminus 
{\cal N}_\epsilon(C) )$ for $i$ large enough.  
Thus by a well-known compactness theorem for 
surfaces with uniformly bounded Gaussian curvature (see, 
for example, \cite{An}), there exists a subsequence of 
the sequence $\{M_i\}_{i=1}^\infty$ which converges in 
$B_{1/\epsilon} \setminus 
{\cal N}_\epsilon(C)$.  The limit of this sequence is nonempty 
if $\epsilon < \delta$, by the 
fourth assumption of the lemma.  Also, by the compactness theorem 
in \cite{An}, if each $M_i$ is embedded, the limit surface is 
embedded.  

By considering a sequence $\{\epsilon_j\}_{j=1}^\infty$ so that 
$\epsilon_j \searrow 0$ as $j \rightarrow \infty$, and by repeatedly 
applying the above argument, we can create a nested sequence 
of convergent subsequences.  The first subsequence 
$\{M_{1i}\}_{i=1}^\infty$ converges in 
$B_{1/{\epsilon_1}} \setminus {\cal N}_{\epsilon_1}(C)$; the second 
subsequence $\{M_{2i}\}_{i=1}^\infty$ is a 
subsequence of $\{M_{1i}\}_{i=1}^\infty$ and converges in 
$B_{1/{\epsilon_2}} \setminus {\cal N}_{\epsilon_2}(C)$; the third 
subsequence $\{M_{3i}\}_{i=1}^\infty$ is a subsequence 
of $\{M_{2i}\}_{i=1}^\infty$ and converges in 
$B_{1/{\epsilon_3}} \setminus {\cal N}_{\epsilon_3}(C)$; and so on.  
By a Cantor diagonalization argument, the subsequence 
$\{M_{ii}\}_{i=1}^\infty$ of the sequence 
$\{M_i\}_{i=1}^\infty$ converges in $\bfR^3$.  
The limit surface $M$ is a surface with boundary $C$. Furthermore, $M$ must 
be stable, for if it were not, it follows that some $M_{ii}$ 
would not be least-area.  
\end{proof}

The following lemma will be used to prove the classification theorem.

\begin{lemma}
Suppose that $M$ is a genus-zero finite-total-curvature complete 
immersed minimal surface.  Suppose that 
$\psi$ is an nontrivial orientation-preserving isometry of $M$.  
Then the set of points and ends of $M$ that are fixed by $\psi$ 
contains at most two elements.  
\end{lemma}

\begin{proof}
The surface $M$ is conformally the sphere with a finite number of 
points removed \cite{Os}; that is, we have a bijective conformal map  
\[ \Phi: \bfC \cup \{\infty\} \setminus \{p_1,...,p_l\} 
\longrightarrow M \; \; \; . \]  The points $\{p_1,...,p_l\}$ 
represent the ends of $M$, and $\Phi$ can be extended 
conformally to the points $\{p_1,...,p_l\}$.
The map $(\Phi)^{-1} \circ \psi \circ \Phi$ extends to a 
bijective conformal map of 
$\bfC \cup \{\infty\}$ to itself.  
If $\psi$ were to fix three or more points or ends of $M$, then 
the extension of $(\Phi)^{-1} \circ \psi \circ \Phi$ would 
have three fixed 
points, and thus would be the identity map.  Therefore $\psi$ 
would be the identity map, a contradiction.
\end{proof}

\section{Previously known minimal surfaces}

In this section, using results from the last section, 
we prove existence of some previously known 
minimal surfaces.  We do this to introduce the 
methods that will be later used to prove existence of the new 
examples, and because nowhere in the literature have 
these old examples explicitly been proven to exist via the conjugate 
surface construction.  

The proofs in \cite{BeRo} depend on the fact that the Jorge-Meeks 
$n$-oids and genus-zero Platonoids 
exist.  The known existence of these genus zero surfaces is used to 
prove the existence of higher genus analogues \cite{BeRo}.  
We prove here the existence of these genus-zero examples.  

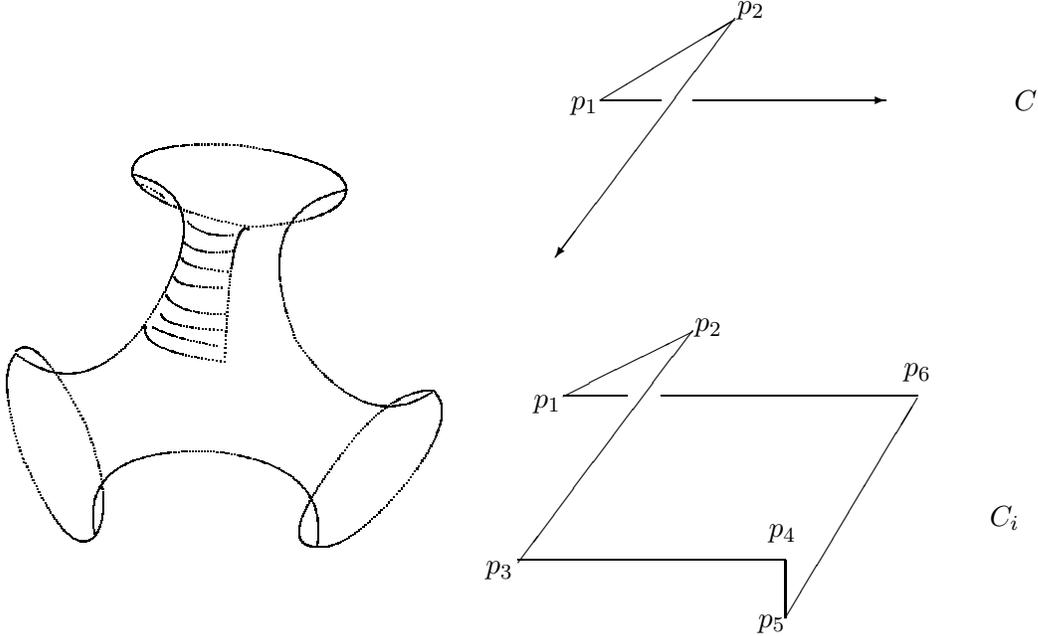
\begin{figure}
\unitlength=1.0pt
\begin{picture}(390.00,232.00)(75.00,535.00)
\put(327.00,733.00){\line(-1,0){23.00}}
\put(304.00,733.00){\line(5,3){51.00}}
\put(314.00,621.00){\line(-1,0){24.00}}
\put(290.00,621.00){\line(2,1){49.00}}
\put(339.00,646.00){\line(-3,-4){66.00}}
\put(273.00,559.00){\line(1,0){101.00}}
\put(374.00,559.00){\line(0,-1){22.00}}
\put(374.00,537.00){\line(3,5){50.00}}
\put(424.00,621.00){\line(-1,0){97.00}}
\put(355.00,764.00){\vector(-3,-4){68.00}}
\put(339.00,733.00){\vector(1,0){73.00}}
\put(298.00,731.00){\makebox(0,0)[cc]{$\small p_1$}}
\put(361.00,767.00){\makebox(0,0)[cc]{$\small p_2$}}
\put(284.00,618.00){\makebox(0,0)[cc]{$\small p_1$}}
\put(345.00,647.00){\makebox(0,0)[cc]{$\small p_2$}}
\put(266.00,555.00){\makebox(0,0)[cc]{$\small p_3$}}
\put(373.00,569.00){\makebox(0,0)[cc]{$\small p_4$}}
\put(369.00,535.00){\makebox(0,0)[cc]{$\small p_5$}}
\put(424.00,630.00){\makebox(0,0)[cc]{$\small p_6$}}
\put(465.00,733.00){\makebox(0,0)[cc]{$\large C$}}
\put(457.00,575.00){\makebox(0,0)[cc]{$\large C_i$}}
\bezier52(127.00,705.00)(131.00,693.00)(171.00,685.00)
\bezier48(171.00,685.00)(207.00,687.00)(208.00,699.00)
\bezier67(208.00,699.00)(208.00,710.00)(167.00,716.00)
\bezier54(167.00,716.00)(126.00,718.00)(127.00,705.00)
\bezier50(83.00,638.00)(91.00,646.00)(109.00,611.00)
\bezier47(109.00,611.00)(121.00,581.00)(113.00,568.00)
\bezier56(113.00,568.00)(103.00,558.00)(86.00,597.00)
\bezier45(86.00,596.00)(75.00,629.00)(83.00,637.00)
\bezier48(241.00,623.00)(251.00,613.00)(229.00,586.00)
\bezier42(229.00,586.00)(208.00,562.00)(197.00,564.00)
\bezier71(197.00,564.00)(178.00,564.00)(211.00,605.00)
\bezier38(211.00,605.00)(227.00,624.00)(241.00,623.00)
\bezier90(208.00,699.00)(170.00,691.00)(189.00,642.00)
\bezier74(189.00,642.00)(215.00,605.00)(241.00,623.00)
\bezier86(113.00,568.00)(106.00,602.00)(158.00,600.00)
\bezier75(158.00,600.00)(202.00,595.00)(197.00,564.00)
\bezier68(131.00,647.00)(108.00,618.00)(83.00,637.00)
\bezier96(131.00,647.00)(164.00,696.00)(127.00,705.00)
\bezier63(171.00,684.00)(164.00,689.00)(162.00,634.00)
\bezier41(162.00,634.00)(130.00,638.00)(132.00,647.00)
\bezier9(131.00,701.00)(139.00,697.00)(139.00,696.00)
\bezier19(148.00,687.00)(153.00,681.00)(165.00,682.00)
\bezier20(147.00,679.00)(152.00,674.00)(165.00,676.00)
\bezier20(146.00,673.00)(146.00,670.00)(163.00,668.00)
\bezier21(143.00,666.00)(144.00,662.00)(161.00,662.00)
\bezier26(140.00,659.00)(144.00,652.00)(162.00,652.00)
\bezier27(138.00,652.00)(139.00,647.00)(161.00,646.00)
\bezier25(135.00,647.00)(142.00,643.00)(159.00,640.00)
\end{picture}

        \caption{A fundamental piece of the Jorge-Meeks surface, boundary 
          contour $C$ of the conjugate of the fundamental piece,
          and the compact boundary contours $C_i$.}
        \label{nnoidproof}
\end{figure}

\begin{theorem}
{\sf (The Jorge-Meeks $n$-oids)} 
For each $n \geq 2$, there exists an 
immersed minimal surface ${\cal J}{\cal M}_0(n)$, 
satisfying the following:
\begin{list}%
{\arabic{num})}{\usecounter{num}\setlength{\rightmargin}{\leftmargin}}
\item ${\cal J}{\cal M}_0(n)$ has genus zero.
\item ${\cal J}{\cal M}_0(n)$ has $n$ catenoid ends with equal weight, 
and the normal vectors at these ends all lie within a plane and are 
symmetrically placed.  
\item The symmetry group of ${\cal J}{\cal M}_0(n)$ is $D_n \times \bfZ_2$.  
\end{list}
\end{theorem}

\begin{proof}
Let $M$ be any immersed smooth minimal surface with a planar 
geodesic $\alpha$ in its boundary, and let Ref($M$) be the 
reflection of $M$ across the plane containing $\alpha$.  
It is well known, by analytic continuation properties 
of minimal surfaces, that $M \cup$ Ref($M$) is a smooth surface 
along $\alpha$, which is now an interior curve of the surface 
\cite{Ka2}.  Therefore the surface ${\cal J}{\cal M}_0(n)$ 
exists if its fundamental piece exists, and its 
fundamental piece would look as in Figure \ref{nnoidproof}.  

The fundamental piece exists if its conjugate surface exists.  
If the conjugate surface exists, it would be a surface with an 
end which is a 90 degree arc of a helicoid end.  The boundary 
$C$ of the conjugate surface, up to a homothety and rigid 
motion of $\bfR^3$, consists of a line segment from $p_1 = (0,0,0)$ 
to $p_2 = (0,\cos(\pi/n),\sin(\pi/n))$, a ray 
pointing in the direction of the positive $x_1$-axis with 
endpoint $p_2$, and a ray 
pointing in the direction of the positive $x_2$-axis with 
endpoint $p_1$.  This follows from the properties of the 
conjugate map, as described in Subsection~2.1.  The points 
$p_1$ and $p_2$ are 
the singular points of the boundary $C$, and the angles that 
$C$ forms at these two points are determined, since the 
conjugate map preserves angles.  

Thus we only need to prove existence of a minimal surface 
with a 90 degree arc of a helicoid end and boundary $C$.  
We do this by finding a sequence 
$\{C_i\}_{i=1}^\infty$ of compact 
contours converging to $C$ and satisfying all the conditions 
of Lemma~2.3 (see Figure \ref{nnoidproof}).  

We now describe the finite contour $C_i$.  
Consider the additional points 
$p_3 = (i,\cos(\pi/n),\sin(\pi/n))$, 
$p_4 = (i,i,\sin(\pi/n))$, $p_5 = (i,i,0)$, 
$p_6 = (0,i,0)$ in $\bfR^3$.  Let $l_i$ be the line segment 
connecting $p_i$ to $p_{i+1}$ for $i = 1, ... ,5$, and let 
$l_6$ be the line segment connecting $p_6$ and $p_1$.  Then 
$C_i = l_1 \cup l_2 \cup l_3 \cup l_4 \cup 
l_5 \cup l_6$.   

The fact that each $C_i$ bounds an 
embedded least-area disk follows from either Theorem 2.1 
or Theorem 2.2.  Thus, by Lemma~2.3, we have a minimal 
surface $M_{\mbox{\footnotesize conj}}$ which is bounded by $C$.  
By Theorem 2.2, each $M_i$ is a connected 
graph over a convex domain $R_i$ in the 
$x_2x_3$-plane, and $R_i \subset R_j$ for $j > i$.  
It follows that $M_{\mbox{\footnotesize conj}}$ is a connected graph, and is therefore 
conformally a disk.  

We do not yet know that the end
of $M_{\mbox{\footnotesize conj}}$ is a 90 degree arc of a helicoid end.  To show this 
we first show that $M_{\mbox{\footnotesize conj}}$ has finite total curvature.  
Choose an orientation on $M_{\mbox{\footnotesize conj}}$, and consider the Gauss 
map $G: M_{\mbox{\footnotesize conj}} \longmapsto S^2$.  Let $P$ be the 
plane containing the points $p_1$, $p_2$ and $p_3$.  
Let Im($M_{\mbox{\footnotesize conj}}$) $\subset S^2$ be the image of
$M_{\mbox{\footnotesize conj}}$ under $G$.  
Note that since $M_{\mbox{\footnotesize conj}}$ is a graph, the image 
Im($M_{\mbox{\footnotesize conj}}$) must lie within a hemisphere.  
Let $N$ be the normal vector to $P$, chosen so that 
$G(p_2) = 
+N$.  Note that $M_{\mbox{\footnotesize conj}}$ lies to one side of $P$ and 
that $C$ makes a 90 degree angle at 
$p_2$.  It follows that $G$ cannot be branched at 
$p_2$.  Furthermore, by comparing the plane $P$ and the surface 
$M_{\mbox{\footnotesize conj}}$ along $C$ and applying the boundary point 
maximum principle, we can conclude that the set 
$G^{-1}(N) \cap C$ consists only of the point $p_2$.  
Thus the branched covering map $G: M_{\mbox{\footnotesize conj}} \longmapsto$ 
Im($M_{\mbox{\footnotesize conj}}$) must be a finite covering map, in fact 
it must have degree one.  In particular, $M_{\mbox{\footnotesize conj}}$ has finite 
total curvature.

Since conjugation is an isometry, we know that the 
fundamental piece $M$ also has finite total curvature.  
We can extend $M$ by reflection to a complete smooth 
finite-total-curvature surface $M_{\mbox{\footnotesize comp}}$.  
Since $M$ is a graph over both the $x_1x_3$-plane and the 
$x_2x_3$-plane (cf. \cite{Ka3}), we can see that the ends of $M_{\mbox{\footnotesize comp}}$ 
are embedded.  Thus they must be of either catenoid or 
planar type \cite{Scn1}.  Suppose they are of planar 
type.  Then $M_{\mbox{\footnotesize conj}}$ must have an end which is 
asymptotic to a plane.  But 
the two boundary rays of $M_{\mbox{\footnotesize conj}}$ do not lie 
in a common plane, so the end of $M_{\mbox{\footnotesize conj}}$ cannot be asymptotic 
to a plane.  Hence the ends of $M_{\mbox{\footnotesize comp}}$ must be of catenoid type.  
This shows that 
the end of $M_{\mbox{\footnotesize conj}}$ is a 90 degree arc of a helicoid end.  

By setting ${\cal J}{\cal M}_0(n) 
= M_{\mbox{\footnotesize comp}}$, the proof is completed.  
\end{proof}

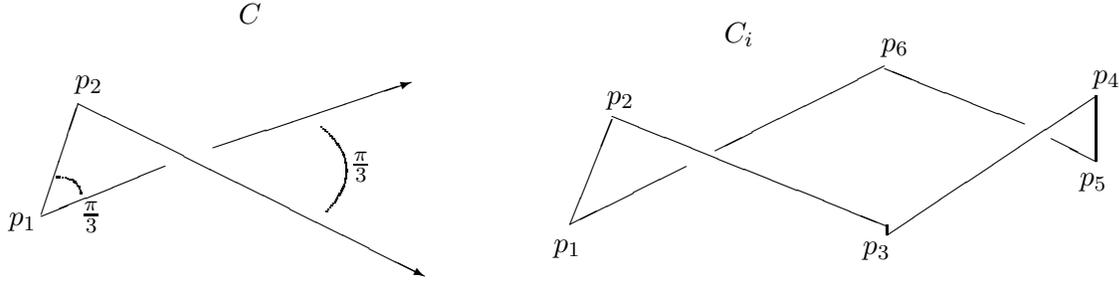
\begin{figure}
\unitlength=1.0pt
\begin{picture}(410.00,116.00)(102.00,644.00)
\put(156.00,685.00){\line(-5,-2){47.00}}
\put(109.00,667.00){\line(1,3){14.00}}
\put(353.00,685.00){\line(-2,-1){44.00}}
\put(309.00,663.00){\line(2,5){16.00}}
\put(325.00,704.00){\line(5,-2){104.00}}
\put(429.00,663.00){\line(0,-1){4.00}}
\put(429.00,659.00){\line(3,2){79.00}}
\put(508.00,712.00){\line(0,-1){25.00}}
\put(508.00,687.00){\line(-2,1){17.00}}
\put(481.00,701.00){\line(-5,2){53.00}}
\put(428.00,723.00){\line(-2,-1){64.00}}
\put(123.00,709.00){\vector(2,-1){131.00}}
\put(174.00,692.00){\vector(3,1){75.00}}
\put(102.00,664.00){\makebox(0,0)[cc]{$\small p_1$}}
\put(127.00,716.00){\makebox(0,0)[cc]{$\small p_2$}}
\put(308.00,654.00){\makebox(0,0)[cc]{$\small p_1$}}
\put(328.00,710.00){\makebox(0,0)[cc]{$\small p_2$}}
\put(425.00,653.00){\makebox(0,0)[cc]{$\small p_3$}}
\put(512.00,718.00){\makebox(0,0)[cc]{$\small p_4$}}
\put(507.00,679.00){\makebox(0,0)[cc]{$\small p_5$}}
\put(432.00,730.00){\makebox(0,0)[cc]{$\small p_6$}}
\put(188.00,743.00){\makebox(0,0)[cc]{$\large C$}}
\put(373.00,735.00){\makebox(0,0)[cc]{$\large C_i$}}
\put(128.00,666.00){\makebox(0,0)[cc]{$\small \frac{\pi}{3}$}}
\put(230.00,685.00){\makebox(0,0)[cc]{$\small \frac{\pi}{3}$}}
\bezier14(115.00,681.00)(123.00,681.00)(124.00,675.00)
\bezier46(215.00,700.00)(233.00,685.00)(217.00,668.00)
\end{picture}

        \caption{The boundary contour $C$ of the conjugate of a 
          fundamental piece of the tetroid, and the 
          compact boundary contours $C_i$.}
        \label{tetproof}
\end{figure}

The same method will be used in all of the 
following proofs (except the proof of Theorem~1.6).  Therefore, 
in the following proofs, we shall only emphasize the differences 
from the previous proof.  We ask the reader to refer to the 
proof of Theorem~3.1 to find information that is left 
unstated in the following arguments.  The proofs of Theorems 
1.1 and 1.2 also include additional period-killing arguments.  

\begin{theorem}
{\sf (The Platonoids)} 
The following genus-zero minimal surfaces with 
\newline catenoid ends exist:
\begin{list}%
{\arabic{num})}{\usecounter{num}\setlength{\rightmargin}{\leftmargin}}

\item {\sf (The genus-zero tetroid)} 
A surface with four ends and symmetry group isomorphic to 
the symmetry group of a tetrahedron.  

\item {\sf (The genus-zero cuboid)} 
A surface with eight ends and symmetry group isomorphic to 
the symmetry group of a cube.  

\item {\sf (The genus-zero octoid)} 
A surface with six ends and symmetry group isomorphic to 
the symmetry group of an octahedron.  

\item {\sf (The genus-zero dodecoid)} 
A surface with twenty ends and symmetry group isomorphic to 
the symmetry group of a dodecahedron.  

\item {\sf (The genus-zero icosoid)} 
A surface with twelve ends and symmetry group isomorphic to 
the symmetry group of an icosahedron.  
\end{list}
\end{theorem}

\begin{proof}
We give here the proof only for the tetroid, as the other four cases 
are similar.  

The surface exists if its fundamental piece exists.  The 
fundamental piece exists if its conjugate surface exists; that is, 
if the noncompact contour $C$ (see 
Figure \ref{tetproof}) bounds a minimal surface with a 60 degree 
arc of a helicoid 
end.  Again, there exists a sequence $\{C_i\}_{i=1}^\infty$ of compact 
contours which bound least-area embedded disks, and which converge to 
$C$ in the sense of Lemma~2.3.  The curve $C_i$ can be chosen to 
be a polygonal contour with vertices $p_1 = (0,0,0)$, 
$p_2 = (-\sqrt{1/8},\sqrt{3/8},1)$, 
$p_3 = p_2 + (i,i/\sqrt{3},0)$, 
$p_4 = (i,i,1)$, $p_5 = (i,i,0)$, and $p_6 = (0,i,0)$, 
connected in the same way as for the proof 
of Theorem 3.1.  Note that $C_i$ makes an angle of 60 degrees at 
$p_1$ and an angle of 90 degrees at $p_2$.  

The result follows just as in the previous proof.  
\end{proof}

\section{Proofs of the main results}

\begin{figure}
\unitlength=1.0pt
\begin{picture}(416.00,181.00)(91.00,563.00)
\put(178.00,609.00){\line(2,3){58.00}}
\put(178.00,603.00){\makebox(0,0)[cc]{$\small p_3$}}
\put(236.00,697.00){\line(4,-5){16.00}}
\put(236.00,701.00){\makebox(0,0)[cc]{$\small p_2$}}
\put(252.00,677.00){\line(0,-1){37.00}}
\put(256.00,683.00){\makebox(0,0)[cc]{$\small p_1$}}
\put(252.00,640.00){\line(3,5){56.00}}
\put(252.00,634.00){\makebox(0,0)[cc]{$\small p_7$}}
\put(308.00,734.00){\line(-3,-1){134.00}}
\put(315.00,734.00){\makebox(0,0)[cc]{$\small p_6$}}
\put(174.00,690.00){\line(0,-1){3.00}}
\put(174.00,687.00){\line(-5,-6){83.00}}
\put(174.00,695.00){\makebox(0,0)[cc]{$\small p_5$}}
\put(91.00,588.00){\line(4,1){88.00}}
\put(91.00,596.00){\makebox(0,0)[cc]{$\small p_4$}}
\put(130.00,655.00){\makebox(0,0)[cc]{$\large C_i$}}
\put(180.00,655.00){\makebox(0,0)[cc]{$\large M_i$}}
\put(330.00,620.00){\makebox(0,0)[cc]{$\large M_{\mbox{\footnotesize conj}}$}}
\put(477.00,665.00){\makebox(0,0)[cc]{$\large C$}}
\put(345.00,563.00){\line(3,5){75.00}}
\put(420.00,688.00){\line(5,-6){18.00}}
\put(420.00,693.00){\makebox(0,0)[cc]{$\small p_2$}}
\put(438.00,666.00){\line(0,-1){38.00}}
\put(444.00,669.00){\makebox(0,0)[cc]{$\small p_1$}}
\put(438.00,628.00){\line(3,5){69.00}}
\put(438.00,621.00){\makebox(0,0)[cc]{$\small p_7$}}

\bezier538(504.00,744.00)(171.00,607.00)(344.00,563.00)
\end{picture}

        \caption{The conjugate of a fundamental piece for the 
          prismoids, and the corresponding finite contours $C_i$.} 
        \label{1.1conj}
\end{figure}
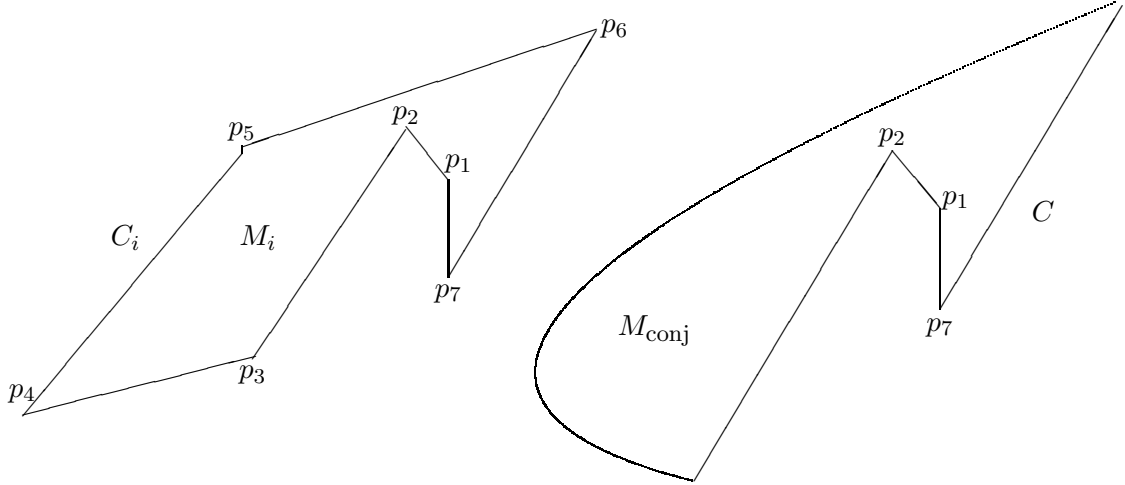

\begin{proof}
({\sf of Theorem 1.1, prismoids})

Using either 
Theorem 2.1 or Theorem 2.2, there exists a sequence of compact 
contours $C_i$ bounding least-area embedded disks $M_i$, 
so that the compact contours converge to the noncompact 
boundary $C$ of the conjugate of a fundamental piece.  
All the conditions of Lemma~2.3 are satisfied.  The Jordan 
contour $C_i$ can be chosen to consist 
of straight line segments from $p_1 = (0,0,0)$ to 
$p_2 = (-t\cos(\pi/n),-t\sin(\pi/n),0)$, then to 
$p_3 = (i,-t\sin(\pi/n),0)$, then to 
$p_4 = (i,-i,t\sin(\pi/n)\cot(\theta)-i\cot(\theta))$, then to 
$p_5 = (-i,-i,-i\cot(\theta)-s)$, then to 
$p_6 = (-i,0,-s)$, then to $p_7 = (0,0,-s)$, and 
then back to $p_1$.  Thus the 
noncompact contour $C$ consists of a line segment from $p_1$ to 
$p_2$, a line segment from $p_1$ to $p_7$, a ray with endpoint 
$p_2$ pointing in the direction of the positive $x_1$-axis, and a 
ray with 
endpoint $p_7$ pointing in the direction of the negative $x_1$-axis.  

Therefore by Lemma 2.3 the conjugate surface $M_{\mbox{\footnotesize conj}}$ of the 
potential fundamental piece exists (see Figure \ref{1.1conj}).  
$M_{\mbox{\footnotesize conj}}$ is conformally a disk, and has an end that is a 180 
degree arc of a helicoid end.  

We now know that the fundamental piece $M$ exists, and 
has a 180 degree arc of a catenoid end.  The boundary 
of this fundamental piece $\partial M$ consists of two finite 
planar geodesics and two infinite planar geodesics.  The two 
infinite planar geodesics lie in parallel planes.  If these 
two infinite planar geodesics lie in the same plane, then 
the entire complete surface exists.  Thus 
there is one period to kill.  The numbers 
$s$, $t > 0$ can be chosen so that the two infinite boundary 
rays of $M_{\mbox{\footnotesize conj}}$ lie in a 
common plane that is perpendicular to the end of 
$M_{\mbox{\footnotesize conj}}$, thus satisfying the conditions 
of Lemma~2.2 (up to a rigid motion).  We conclude by Lemma 2.2 that, 
for these values of $s$ and $t$, 
$M$ extends to a complete finite-total-curvature minimal 
surface.  Thus the period is zero, and the proof is completed.  
\end{proof}

\begin{remark}
({\sf Jorge-Meeks $n$-oid fence})

By the method of the first two paragraphs in the 
previous proof, we can construct embedded 
minimal disks $M_s, \; s \in (0,\infty)$ 
with the following properties:  $M_s$ is bounded by a 
straight line segment from $(0,0,0)$ to $(0,-\cos(\pi/n),
-\sin(\pi/n))$, a straight line segment from $(0,0,0)$ to 
$(-s,0,0)$, a ray with the endpoint $(0,-\cos(\pi/n),-\sin(\pi/n))$ 
pointing in the direction of the 
negative $x_1$-axis, and a ray with the endpoint $(-s,0,0)$ pointing 
in the direction of the negative $x_2$-axis.  Furthermore, 
$M_s$ has a single end that is 
90 degree arc of a helicoid end, and $M_s$ is a graph over the 
region $\{(0,x_2,x_3) \; | \; \, \, -\sin(\pi/n) < x_3 < 0, 
\, x_2 < \cot(\pi/n) \cdot x_3 \}$.

Consider the conjugate surface to $M_s$ and this conjugate 
surface's extension by reflection across boundary planar geodesics 
to a complete minimal surface.  We call the 
resulting surface the Jorge-Meeks $n$-oid fence (see Figure 
\ref{between} (2)).  
It is a surface with translational symmetry in one direction.  
The portion of the surface which generates the entire surface 
under the translation has $n$ symmetrically placed ends, just 
as for the Jorge-Meeks $n$-oid.  The complete surface looks
like an infinite collection of $n$-oids regularly spaced along 
a single direction, with each pair of adjacent ``$n$-oids'' connected 
by a handle.

There is a one-parameter family of Jorge-Meeks $n$-oid fences, 
one surface for each value of $s > 0$.  In the case $n=2$ we have
the catenoid fence (cf. \cite{Ka3}).
\end{remark}

\begin{figure}
\hspace{-1.25in}
\unitlength=1.0pt
\begin{picture}(418.00,257.00)(6.00,540.00)
\put(420.00,797.00){\line(-2,-3){60.00}}
\put(360.00,707.00){\line(5,3){59.00}}
\put(419.00,743.00){\line(-2,-3){42.00}}
\put(377.00,680.00){\line(0,-1){50.00}}
\put(377.00,630.00){\line(-2,-3){60.00}}
\put(425.00,741.00){\makebox(0,0)[cc]{$\small p_1$}}
\put(383.00,675.00){\makebox(0,0)[cc]{$\small p_2$}}
\put(384.00,633.00){\makebox(0,0)[cc]{$\small p_3$}}
\put(372.00,581.00){\makebox(0,0)[cc]{$\large C$}}
\put(239.00,623.00){\makebox(0,0)[cc]{$\large M_{\mbox{\footnotesize conj}}$}}
\put(354.00,707.00){\makebox(0,0)[cc]{$\small p_8$}}
\bezier799(318.00,540.00)(6.00,540.00)(420.00,797.00)
\end{picture}
        \caption{The conjugate of a fundamental piece for the 
          higher genus prismoids.}
        \label{1.2conj}
\end{figure}
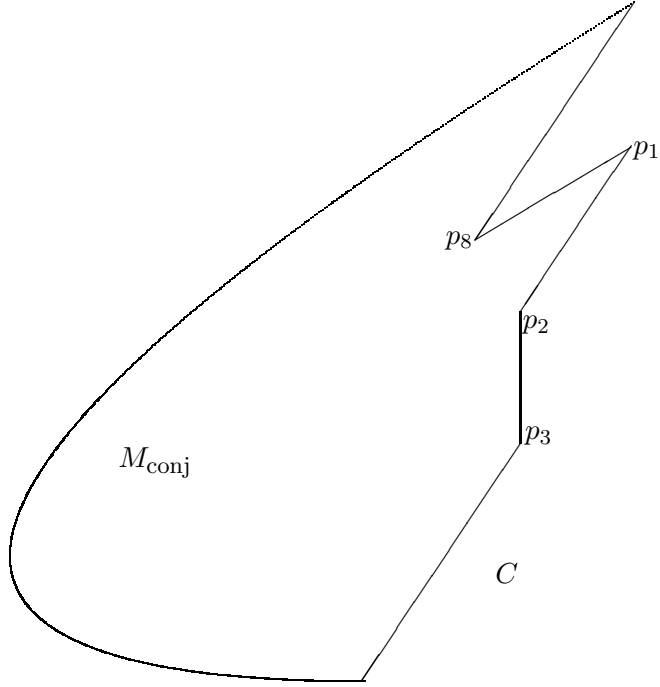

\begin{proof}
({\sf of Theorem 1.2, higher genus prismoids})

Again the conjugate surface of the fundamental piece 
exists by Theorem 2.1 or Theorem 2.2, followed 
by Lemma 2.3.  The finite polygonal contours $C_i$ 
can be chosen to consist of straight line segments from 
$p_1 = (-s,0,0)$ to $p_2 = (0,0,0)$, then to $p_3 = (0,0,-t)$, 
then to $p_4 = (i,0,-t)$, then to 
$p_5 = (i,-i,-t-i\cot(\theta))$, then to 
$p_6 = (-i,-i,-i\cot(\theta)+u\sin(\pi/n)\cot(\theta))$, 
then to 
$p_7 = (-i,-u\sin(\pi/n),0)$, then to 
$p_8 = (u\cos(\pi/n)-s,-u\sin(\pi/n),0)$, 
and then back to $p_1$.  
Thus the limit contour $C$ consists of a line segment from 
$p_1$ to $p_8$, a line segment from $p_1$ to $p_2$, a line 
segment from $p_2$ to $p_3$, a ray with endpoint $p_8$ pointing 
in the direction of the negative $x_1$-axis, and a ray with 
endpoint $p_3$ pointing in the direction of the positive 
$x_1$-axis.  Furthermore, the limit surface $M_{\mbox{\footnotesize conj}}$ 
bounded by $C$ has a normal vector at its end which makes 
an angle of $\theta$ with a horizontal plane (see Figure 
\ref{1.2conj}).  The 
conjugate surface $M_{\mbox{\footnotesize conj}}$ is conformally a disk, and has 
an end that is a 180 degree arc of a helicoid end.  

Here we have two periods to kill.   Let $\alpha_1$ be the 
unbounded boundary curve on the fundamental piece $M$ that corresponds 
(via the conjugate map) to the boundary ray of $M_{\mbox{\footnotesize conj}}$ with 
endpoint $p_8$.  Let $\alpha_2$ be the 
unbounded boundary curve of $M$ corresponding to the boundary ray 
of $M_{\mbox{\footnotesize conj}}$ with endpoint $p_3$.  Let $\alpha_3$ be the bounded 
boundary curve of $M$ corresponding to the boundary line 
segment of $M_{\mbox{\footnotesize conj}}$ with endpoints $p_1$ and $p_2$.  To kill 
both periods we must show there exist choices of $s,t,u > 0$ 
so that $\alpha_1$, $\alpha_2$, and $\alpha_3$ all lie within 
a common plane.

Values can be chosen for $u,t$ so that the two boundary rays of 
$M_{\mbox{\footnotesize conj}}$ lie in a common plane that is perpendicular to the 
end of $M_{\mbox{\footnotesize conj}}$.  Then, by Lemma 2.2, we conclude that $\alpha_1$ and 
$\alpha_2$ lie in a common plane $P$.  Note that these values of $u,t$ are 
independent of the value of $s$.  

We now show that for some value of $s > 0$, the curve $\alpha_3$ 
also lies in $P$.  The surface $M_{\mbox{\footnotesize conj}}$ varies smoothly in 
$s \in [0,\infty)$, and is a graph over a fixed region in the 
$x_2x_3$-plane for all $s \in [0,\infty)$.  
As $s \rightarrow 0$, $M_{\mbox{\footnotesize conj}}$ converges 
smoothly to the conjugate surface of a fundamental piece of 
${\cal P}_0(2n,\theta)$.  We shall refer to this fundamental piece of 
${\cal P}_0(2n,\theta)$ as $FM$.  

\begin{figure}
\unitlength=1.0pt
\begin{picture}(398.00,84.00)(114.00,679.00)
\put(114.00,707.00){\line(1,0){398.00}}
\put(142.00,707.00){\circle*{4}}
\put(166.00,742.00){\circle*{4}}
\put(337.00,707.00){\circle*{4}}
\put(244.00,736.00){\circle*{4}}
\put(405.00,707.00){\circle*{4}}
\put(482.00,707.00){\circle*{4}}
\put(157.00,686.00){\makebox(0,0)[cc]{$s \approx 0$}}
\put(291.00,735.00){\makebox(0,0)[cc]{\small half-circle}}
\put(292.00,686.00){\makebox(0,0)[cc]{$s >> 0$}}
\put(444.00,686.00){\makebox(0,0)[cc]{period killed for some $s > 0$}}
\bezier61(142.00,707.00)(144.00,750.00)(162.00,750.00)
\bezier10(162.00,750.00)(167.00,745.00)(166.00,742.00)
\bezier80(405.00,707.00)(405.00,747.00)(445.00,750.00)
\bezier80(445.00,750.00)(484.00,748.00)(482.00,707.00)
\bezier82(337.00,707.00)(336.00,748.00)(295.00,750.00)
\bezier84(295.00,750.00)(250.00,746.00)(253.00,707.00)
\bezier111(337.00,707.00)(344.00,761.00)(287.00,763.00)
\bezier65(287.00,763.00)(246.00,760.00)(244.00,736.00)
\end{picture}

        \caption{The conjugate 
          images of the line segment from $p_2$ to $p_3$ 
          for various values of $s$, in the proof of Theorem 
          1.2.}  
        \label{compare}
\end{figure}
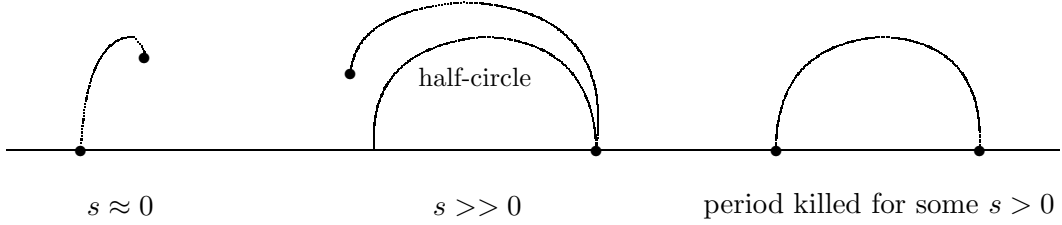

To describe the behavior of $M_{\mbox{\footnotesize conj}}$ as 
$s \rightarrow \infty$, we first describe a portion of a 
helicoid.  Let $H$ be a portion of a helicoid with a 
single end that traverses 180 degrees, and has 
boundary consisting of the line segment from $p_2$ to $p_3$, the 
ray with endpoint $p_3$ pointing the direction of the 
positive $x_1$-axis, and the ray with endpoint $p_2$ pointing the 
direction of the negative $x_1$-axis.  Choose $H$ so that it 
is a nonempty graph over $\{(x_1,x_2,0) \in \bfR^3 \; | \; x_2 < 0\}$, 
thus $H$ is unique.  
As $s \rightarrow \infty$, ${\cal P}_0(2n,\theta)$ converges smoothly to 
a surface which looks similar to $H$, except that its end is 
``slanted'' by the angle $\theta$, hence we call this limit 
surface $SH$.  Note that $\partial SH = \partial H$.  

By the maximum principle, $H$ and $SH$ are disjoint in their 
interiors, and $H$ lies above $SH$.  (This can be argued 
rigorously by a ``sliding'' argument, see \cite{BeRo}.)  
Thus, as one travels downward along the boundary 
line from $p_2$ to $p_3$, the normal vector of $H$ must be 
turning ahead of the normal vector of $SH$.  The same is 
then true of the corresponding boundary curves of the 
conjugate surfaces (with respect to arc length along these two 
planar geodesics).  Since the catenoid is the conjugate surface 
of the helicoid, the conjugate image of the boundary
line from $p_2$ to $p_3$ with respect to $H$ is a half-circle.  
The conjugate image of the boundary 
line from $p_2$ to $p_3$ with respect to $SH$ is not a half-circle, 
but it has the same length as the half-circle, since conjugation 
is an isometry.

It follows that for large values of $s$, the curve $\alpha_3$ 
lies strictly to one particular side of $P$.  By an examination of the 
placement of $FM$ in $\bfR^3$, we see that for values of $s$ 
close to zero, the curve $\alpha_3$
lies strictly to the other side of $P$.  Thus, by the 
intermediate-value-theorem, there exists some value of $s$ 
so that $\alpha_3 \subset P$ (see Figure \ref{compare}).  
\end{proof}

\begin{figure}
\unitlength=1.0pt
\begin{picture}(423.00,245.00)(91.00,535.00)
\put(100.00,720.00){\line(1,2){30.00}}
\put(100.00,720.00){\line(0,-1){160.00}}
\put(100.00,560.00){\line(1,2){30.00}}
\put(260.00,720.00){\line(0,-1){180.00}}
\put(100.00,720.00){\line(1,0){160.00}}
\put(275.00,570.00){\line(0,1){180.00}}
\put(275.00,750.00){\line(-1,0){160.00}}
\put(264.00,620.00){\line(1,0){8.00}}
\put(279.00,620.00){\line(1,0){162.00}}
\put(260.00,540.00){\line(1,2){39.00}}
\put(302.00,624.00){\line(1,2){58.00}}
\put(440.00,620.00){\line(1,2){70.00}}
\put(510.00,760.00){\line(0,-1){20.00}}
\put(510.00,740.00){\line(-1,0){7.00}}
\put(360.00,740.00){\line(1,0){135.00}}
\put(495.00,740.00){\line(-1,0){1.00}}
\put(130.00,780.00){\line(0,-1){24.00}}
\put(130.00,745.00){\line(0,-1){19.00}}
\put(100.00,650.00){\line(1,2){30.00}}
\put(101.00,650.00){\line(3,-2){159.00}}
\put(130.00,717.00){\line(0,-1){81.00}}
\put(130.00,620.00){\line(0,1){8.00}}
\put(130.00,620.00){\line(1,0){11.00}}
\put(152.00,620.00){\line(1,0){101.00}}
\put(260.00,720.00){\line(1,2){15.00}}
\put(290.00,620.00){\line(0,-1){20.00}}
\put(297.00,597.00){\makebox(0,0)[cc]{{$\small p_1$}}}
\put(287.00,629.00){\makebox(0,0)[cc]{{$\small p_2$}}}
\put(447.00,614.00){\makebox(0,0)[cc]{{$\small p_3$}}}
\put(504.00,771.00){\makebox(0,0)[cc]{{$\small p_4$}}}
\put(514.00,734.00){\makebox(0,0)[cc]{{$\small p_5$}}}
\put(353.00,751.00){\makebox(0,0)[cc]{{$\small p_6$}}}
\put(283.00,567.00){\makebox(0,0)[cc]{{$\small p_7$}}}
\put(264.00,535.00){\makebox(0,0)[cc]{{$\small p_8$}}}
\put(97.00,553.00){\makebox(0,0)[cc]{{$\small p_9$}}}
\put(91.00,726.00){\makebox(0,0)[cc]{{$\small p_{10}$}}}
\put(106.00,758.00){\makebox(0,0)[cc]{{$\small p_{11}$}}}
\put(140.00,780.00){\makebox(0,0)[cc]{{$\small p_{12}$}}}
\put(122.00,624.00){\makebox(0,0)[cc]{{$\small p_{13}$}}}
\put(248.00,711.00){\makebox(0,0)[cc]{{$\small p_{14}$}}}
\put(281.00,756.00){\makebox(0,0)[cc]{{$\small p_{15}$}}}
\put(139.00,710.00){\makebox(0,0)[cc]{{$\small p_{16}$}}}
\put(91.00,654.00){\makebox(0,0)[cc]{{$\small p_{17}$}}}
\bezier66(100.00,560.00)(133.00,559.00)(165.00,550.00)
\bezier98(165.00,550.00)(182.00,539.00)(260.00,540.00)
\end{picture}

        \caption{The construction of $\hat{M}_i$ in the 
          proof of the genus-zero $n$-oids plus two ends.}
        \label{1.3conj}
\end{figure}

\begin{proof}
({\sf of Theorem 1.3, genus-zero $n$-oids plus two ends})

For this proof there is no period problem, but since the conjugate of 
the fundamental piece is not a graph over a {\em convex} plane domain,
Nitsche's theorem cannot be applied to show existence of the conjugate 
piece.  Thus Theorem 2.1 must be used for this, followed by Lemma 2.3.  

We describe now the 
construction of compact 3-manifolds $\hat{M}_i \subset \bfR^3$ 
and the finite polygonal contours $C_i \subset \partial 
\hat{M}_i$.  We construct $\hat{M}_i$ so that it satisfies the hypotheses 
of Theorem 2.1, and thus $C_i$ bounds a least-area 
embedded disk $M_i$.  The result follows as in the previous proofs.  

Assume that $w \geq n/2$, thus 
$w/{2n} \geq 1/4$.  Later we shall need to assume 
that $w > c(n) \geq n/2$, for some constant $c(n)$ 
depending only on $n$.  

The skeletal structure of $\hat{M}_i$ is given in Figure \ref{1.3conj}.  
Let $p_1 = (0,0,0)$, $p_2 = (0,0,1/4)$, $p_3 = (0,i,1/4)$, 
$p_4 = (-i+1/4,i,1/4)$, 
$p_5 = (-i+1/4,i,0)$, 
$p_6 = (-i+1/4,0,0)$, 
$p_{7} = (1/4,0,0)$, 
$p_{8} = (w/{2n},0,0)$, 
$p_{9} = (w/{2n},-i,1/4)$, 
$p_{10} = (w/{2n},-i,i)$, 
$p_{11} = (1/4,-i,i)$, $p_{12} = (0,-i,i)$, $p_{13} = 
(0,-i,1/4)$, $p_{14} = (w/{2n},0,i)$, $p_{15} = 
(1/4,0,i)$, 
$p_{16} = (0,-i,i\tan(\pi/n))$, and 
$p_{17} = (w/{2n},-i,i\tan(\pi/n))$.  

Consider the polygonal contour defined by connecting the 
following vertices by line segments: $p_2$ to $p_7$, 
$p_7$ to $p_{6}$ (through $p_1$), $p_{6}$ to $p_{5}$, $p_{5}$ to 
$p_4$, $p_4$ to $p_3$, and $p_3$ back to $p_2$.  By Theorem 2.1 
or Theorem 2.2, this contour bounds a minimal graph $M^\prime$.  
Let Rot($M^\prime$) be the surface that results from rotating $M^\prime$ 
by 180 degrees about the line through $p_2$ and $p_7$.  Then 
$M^\prime \cup$Rot($M^\prime$) is a smooth disk, which we shall call 
$M_{i1}$.  Let $C_{i1}$ be the boundary of $M_{i1}$.  
Thus $C_{i1}$ is the polygonal contour defined by connecting the 
following vertices by line segments: $p_3$ to $p_4$, 
$p_4$ to $p_{5}$, $p_5$ to $p_{6}$, $p_{6}$ to 
$p_7$ (through $p_1$), $p_7$ to $p_{15}$, $p_{15}$ to $p_{11}$, 
$p_{11}$ to $p_{12}$, $p_{12}$ to $p_{13}$ (through $p_{16}$), 
and $p_{13}$ back to $p_3$ (through $p_2$).  

Consider the polygonal contour defined by connecting the 
following vertices: $p_2$ to $p_1$, 
$p_1$ to $p_{6}$, $p_{6}$ to $p_{5}$, $p_{5}$ to 
$p_4$, $p_4$ to $p_3$, and $p_3$ back to $p_2$.  Again by 
either Theorem 2.1 
or Theorem 2.2, this contour bounds a minimal graph $M^\prime$.  
Let Rot($M^\prime$) be the surface that results from rotating $M^\prime$ 
by 180 degrees about the line through $p_2$ and $p_1$.  Thus we 
have the smooth disk $M^\prime \cup$Rot($M^\prime$).  Let $P$ be the plane 
$\{x_1 = w/{2n}\}$.  Then $( \, M^\prime \cup$Rot$(M^\prime) \, ) 
\setminus P$ has two components.  Consider the component in the
region $\{x_1 \leq w/{2n}\}$, and name it $M_{i2}$.  
Let $C_{i2}$ be the boundary of $M_{i2}$.  
Thus $C_{i2}$ is the polygonal contour defined by connecting 
$p_3$ to $p_4$, $p_4$ to $p_{5}$, $p_5$ to $p_{6}$, $p_{6}$ to 
$p_8$ (through $p_1$ and $p_7$), $p_8$ to $p_{9}$, $p_{9}$ to $p_{13}$, 
and $p_{13}$ back to $p_{3}$ (through $p_{2}$).  All these connections
are made by straight lines, except for the connection from $p_8$ to 
$p_9$, which is made by the nonstraight planar curve 
$\alpha = P \cap$Rot($M^\prime$).  

For our construction to be valid, we need that the line segment 
from $p_8$ to $p_{17}$ lies above $\alpha$, with respect to the 
positive $x_3$ direction.  Note that 
the surface $M^\prime \cup$Rot($M^\prime$) in the previous paragraph 
converges to a portion of a helicoid as $i \rightarrow \infty$.  
Thus the function ${dx_3}/{dx_1}$ along $\alpha$ approaches 
zero as $w$ and $i$ become large.  It follows that 
there exists a positive number $c(n)$ so that if $w > c(n)$, then 
the line segment from $p_8$ to $p_{17}$ lies above $\alpha$.  

Let $M_{i3}$ be the plane rectangle with vertices $p_7$, $p_8$, 
$p_{14}$, and $p_{15}$.  
Let $M_{i4}$ be the plane rectangle with vertices $p_{9}$, $p_{10}$, 
$p_{12}$, and $p_{13}$.  
Let $M_{i5}$ be the plane rectangle with vertices $p_{10}$, $p_{11}$, 
$p_{14}$, and $p_{15}$.  
Let $M_{i6}$ be the plane region in $P$ bounded by the line 
segment from $p_8$ to 
$p_{14}$, the line segment from $p_{14}$ to $p_{10}$, the line 
segment from $p_{10}$ to 
$p_9$, and the curve $\alpha$ from $p_8$ to $p_9$.  

By the maximum principle, the surfaces $M_{ij}$, $j = 1,...,6$, 
are pairwise 
disjoint in their interiors.  Their union forms the boundary of 
a compact 3-manifold $\hat{M}_i$, and it is clear that $\hat{M}_i$ 
satisfies all the conditions of Theorem 2.1.  Let $C_i$ be the 
polygonal Jordan curve connecting the following vertices by 
line segments: $p_3$ to $p_4$, $p_4$ to $p_{5}$, $p_5$ to $p_{6}$, 
$p_{6}$ to $p_8$ (through $p_1$ and $p_7$), $p_8$ to $p_{17}$, 
$p_{17}$ to $p_{16}$, $p_{16}$ to $p_{13}$, and $p_{13}$ back to 
$p_{3}$ (through $p_{2}$).  Since $C_i$ is a curve in the boundary 
of $\hat{M}_i$, we have by Theorem 2.1 that $C_i$ bounds a 
least-area embedded disk, as desired.  
\end{proof}

\begin{proof}
({\sf of Theorem 1.4, 2$n$-oids with alternating weights})

If $w < 1$, we can apply a homothety with dilation factor 
$1/w$ to get an equivalent surface, but with $w > 1$.  
So we may assume $w \geq 1$.  And since 
the case $w=1$ is the known Jorge-Meeks $2n$-oid, we may further 
assume $w > 1$.  

There is no period problem here, but Nitsche's theorem again 
does not apply, so Theorem 2.1 followed by Lemma 2.3 is necessary.  
We describe now the 
construction of compact 3-manifolds $\hat{M}_i \subset \bfR^3$ 
and the finite polygonal contours $C_i \subset \partial 
\hat{M}_i$.  We show that $\hat{M}_i$ satisfies the hypotheses 
of Theorem 2.1, and thus $C_i$ bounds a least-area 
embedded disk $M_i$, and the result follows as before.  

The skeletal structure of $\hat{M}_i$ is given in Figure \ref{1.4}.  
Let $p_2 = (0,0,0)$, $p_3 = (-i,0,0)$, $p_4 = (-i,i,0)$, 
$p_5 = (-i,i,1/4)$, $p_6 = (-i,i,w/4)$, 
$p_7 = (0,i,w/4)$, $p_{12} = (0,i,1/4)$, 
$p_{16} = (i,0,0)$, $p_9 = (0,-\cot(\pi/{2n})/4,
1/4)$, 
$p_1 = p_9 + s(0,\cos(\pi/n),-\sin(\pi/n))$ with $s = 
(4\sin(\pi/{2n})\cos(\pi/{2n}))^{-1}$, 
$p_8 = p_9 + (1-w)(4\sin(\pi/n))^{-1}
(0,\cos(\pi/n),-\sin(\pi/n))$, 
$p_{13} = p_9 + (i + (1/4)\cot(\pi/{2n}))
(0,\cos(\pi/n),-\sin(\pi/n))$, 
$p_{14} = p_{13} + (i,0,0)$, $p_{15} = 
p_{16} + i(0,\cos(\pi/n),-\sin(\pi/n))$, 
$p_{10} = p_9 + (-i,0,0)$, and $p_{11} = p_8 + (-i,0,0)$.  

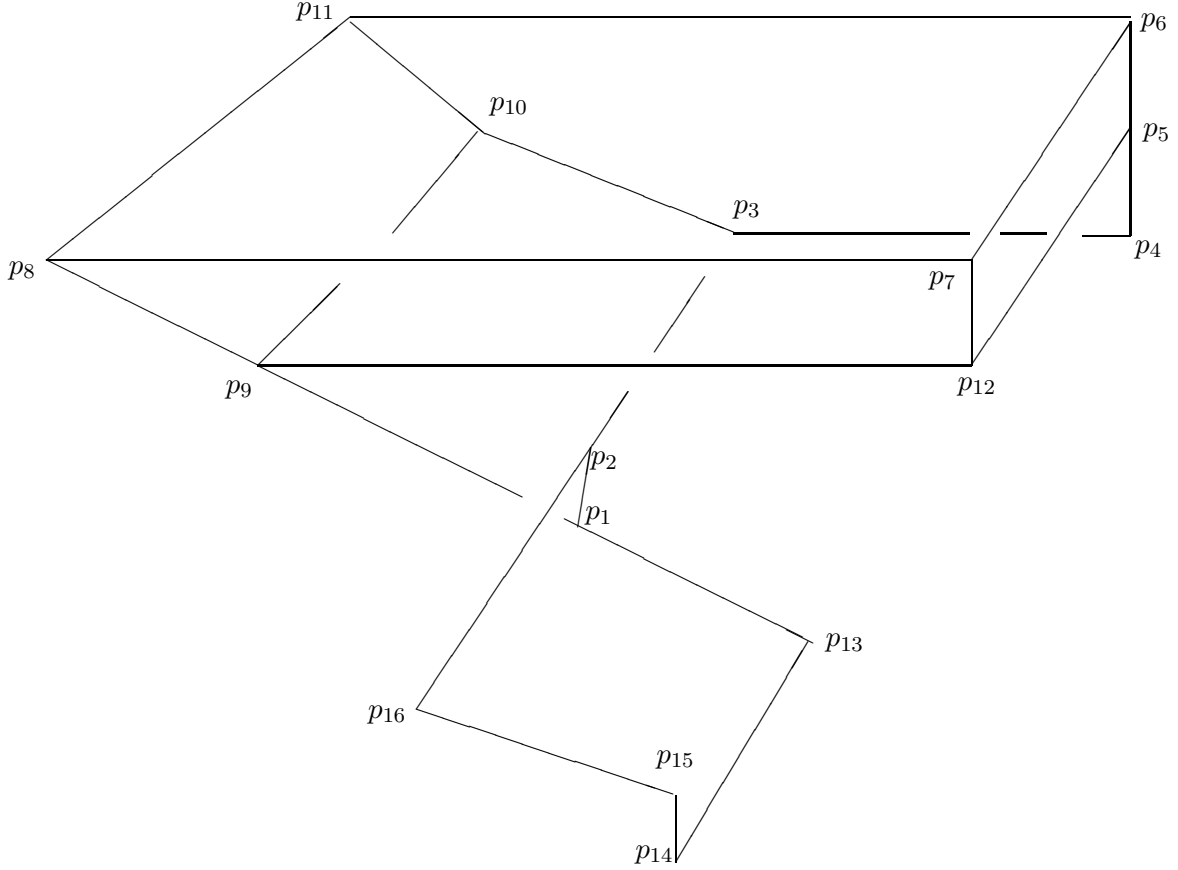
\begin{figure}
\unitlength=1.0pt
\begin{picture}(429.00,322.00)(81.00,492.00)
\put(90.00,720.00){\line(1,0){350.00}}
\put(440.00,720.00){\line(2,3){60.00}}
\put(90.00,720.00){\line(2,-1){180.00}}
\put(230.00,550.00){\line(2,3){80.00}}
\put(170.00,680.00){\line(1,0){270.00}}
\put(440.00,680.00){\line(0,1){1.00}}
\put(440.00,681.00){\line(0,1){39.00}}
\put(440.00,680.00){\line(2,3){60.00}}
\put(500.00,810.00){\line(0,-1){80.00}}
\put(500.00,730.00){\line(0,-1){1.00}}
\put(500.00,729.00){\line(-1,0){18.00}}
\put(468.00,730.00){\line(-1,0){17.00}}
\put(350.00,730.00){\line(1,0){89.00}}
\put(320.00,685.00){\line(2,3){19.00}}
\put(90.00,720.00){\line(5,4){115.00}}
\put(205.00,812.00){\line(1,0){295.00}}
\put(286.00,622.00){\line(2,-1){94.00}}
\put(230.00,550.00){\line(3,-1){97.00}}
\put(170.00,680.00){\line(1,1){31.00}}
\put(328.00,517.00){\line(0,-1){25.00}}
\put(328.00,492.00){\line(3,5){50.00}}
\put(291.00,619.00){\line(1,6){5.00}}
\put(205.00,810.00){\line(6,-5){51.00}}
\put(256.00,768.00){\line(5,-2){94.00}}
\put(221.00,730.00){\line(5,6){32.00}}
\put(348.00,726.00){\line(5,6){3.00}}
\put(299.00,623.00){\makebox(0,0)[cc]{$\small p_1$}}
\put(301.00,644.00){\makebox(0,0)[cc]{$\small p_2$}}
\put(355.00,739.00){\makebox(0,0)[cc]{$\small p_3$}}
\put(507.00,724.00){\makebox(0,0)[cc]{$\small p_4$}}
\put(510.00,768.00){\makebox(0,0)[cc]{$\small p_5$}}
\put(509.00,810.00){\makebox(0,0)[cc]{$\small p_6$}}
\put(192.00,814.00){\makebox(0,0)[cc]{$\small p_{11}$}}
\put(81.00,716.00){\makebox(0,0)[cc]{$\small p_8$}}
\put(163.00,671.00){\makebox(0,0)[cc]{$\small p_9$}}
\put(265.00,778.00){\makebox(0,0)[cc]{$\small p_{10}$}}
\put(442.00,672.00){\makebox(0,0)[cc]{$\small p_{12}$}}
\put(392.00,575.00){\makebox(0,0)[cc]{$\small p_{13}$}}
\put(320.00,495.00){\makebox(0,0)[cc]{$\small p_{14}$}}
\put(328.00,531.00){\makebox(0,0)[cc]{$\small p_{15}$}}
\put(219.00,548.00){\makebox(0,0)[cc]{$\small p_{16}$}}
\put(429.00,712.00){\makebox(0,0)[cc]{$\small p_7$}}
\end{picture}
        \caption{The construction of $\hat{M}_i$ in 
          the proof of the 2$n$-oids with alternating weights.} 
        \label{1.4}
\end{figure}

Let $C_{i1}$ be the polygonal contour defined by connecting the 
following vertices by line segments: $p_8$ to $p_9$, 
$p_9$ to $p_{10}$, $p_{10}$ to $p_{11}$, and $p_{11}$ 
back to $p_8$.  Let $C_{i2}$ be the contour connecting 
$p_3$ to $p_4$, $p_4$ to $p_6$ (through $p_5$), $p_6$ 
to $p_{11}$, $p_{11}$ to $p_{10}$, and $p_{10}$ back to $p_3$.  
Let $C_{i3}$ be the contour connecting 
$p_9$ to $p_{12}$, $p_{12}$ to $p_7$, $p_7$ 
to $p_8$, and $p_{8}$ back to $p_9$.  
Let $C_{i4}$ be the contour connecting 
$p_{12}$ to $p_5$, $p_5$ to $p_6$, $p_6$ 
to $p_7$, and $p_7$ back to $p_{12}$.  
Let $C_{i5}$ be the contour connecting 
$p_7$ to $p_{6}$, $p_{6}$ to $p_{11}$, $p_{11}$ 
to $p_8$, and $p_{8}$ back to $p_7$.  
These five contours are all bound plane regions, 
which we shall call $M_{i1}, M_{i2}, M_{i3}, M_{i4}$, and $M_{i5}$, 
respectively.  

Let $C_{i6}$ be the contour connecting 
$p_1$ to $p_2$, $p_2$ to $p_3$, $p_3$ to $p_{10}$, $p_{10}$ 
to $p_9$, and $p_9$ back to $p_1$.  
Let $C_{i7}$ be the contour connecting 
$p_1$ to $p_2$, $p_2$ to $p_{16}$, 
$p_{16}$ to $p_{15}$, $p_{15}$ to $p_{14}$, $p_{14}$ 
to $p_{13}$, and $p_{13}$ back to $p_1$.  
By Theorem 2.1 or Theorem 2.2, these two contours bound minimal 
graphs $M_{i6}, M_{i7}$, respectively.  

Consider the contour connecting $p_2$ to $p_3$, $p_3$ to $p_4$, $p_4$ 
to $p_5$, $p_5$ to $p_{12}$, $p_{12}$ to $p_9$, and $p_9$ back to 
$p_2$.  By Theorem 2.1 or Theorem 2.2, this contour bounds a minimal 
graph $M^\prime$.  We rotate $M^\prime$ by 180 degrees 
about the line through $p_2$ and $p_9$ to obtain the surface 
Rot$(M^\prime)$.  Then $M_{i8} = M^\prime \cup$Rot$(M^\prime)$ is 
a smooth embedded minimal disk.  Let $C_{i8}$ be the contour connecting
$p_3$ to $p_4$, $p_4$ to $p_5$, $p_5$ to $p_{12}$, $p_{12}$ 
to $p_9$, $p_9$ to $p_{13}$ (through $p_1$), 
$p_{13}$ to $p_{14}$, $p_{14}$ to $p_{15}$, 
$p_{15}$ to $p_{16}$, and $p_{16}$ back to $p_3$ (through $p_2$).  
Then $\partial M_{i8} = C_{i8}$.

The surfaces $M_{ij}, j = 1, ... ,8$ are pairwise disjoint 
in their interiors, by the maximum principle.  Therefore, the 
union of these eight surfaces forms the boundary of a 
region $\hat{M}_i$ in $\bfR^3$.  Let $C_i$ be the contour connecting 
$p_3$ to $p_4$, $p_4$ to $p_6$ (through $p_5$), $p_6$ to $p_7$, $p_7$ 
to $p_8$, $p_8$ to $p_{13}$ (through $p_9$ and $p_1$), 
$p_{13}$ to $p_{14}$, $p_{14}$ to $p_{15}$, 
$p_{15}$ to $p_{16}$, and $p_{16}$ back to $p_3$ (through $p_2$).  
We wish to show that $\hat{M}_i$ satisfies the conditions of 
Theorem 2.1, thus $C_i$ would bound an embedded least-area disk 
$M_i$.  
$\hat{M}_i$ clearly satisfies the conditions of Theorem 2.1, except 
possibly along the line segment between $p_1$ and $p_2$, and 
possibly along the line segment between $p_2$ and $p_3$.  

Let $l_2$ be the line through $p_2$ and $p_3$.  
We now start to rotate $M_{i6}$ about $l_2$ (in the 
clockwise direction with respect to the vector from $p_2$ 
to $p_3$).  
By the maximum principle, the first moment of contact between 
$M_{i8}$ and the interior of $M_{i6}$ cannot occur as a 
tangential contact along $l_2$ and cannot occur at a point 
in the interior of $M_{i8}$.  Thus it must occur as a 
nontangential contact along the line segment from $p_{12}$ to $p_9$, 
which occurs only after the rotation 
has traversed an arc of 180 degrees.  This implies 
that the angles between $M_{i6}$ and $M_{i8}$ 
along $l_2$ must be at most 180 degrees 
with respect to the interior of $\hat{M}_i$.  Thus the conditions 
of Theorem 2.1 are satisfied along the line segment from 
$p_2$ to $p_3$.  

Let $l_1$ be the line through $p_1$ and $p_2$.  Note that 
$l_1$ is perpendicular to the line through $p_2$ and 
$p_9$.  
Here we start to rotate $M_{i6}$ about $l_1$ (in the 
clockwise direction with respect to the vector from $p_1$ 
to $p_2$).  By arguing just as in the previous paragraph, 
we conclude that the angles between $M_{i6}$ and $M_{i7}$ 
along $l_1$ are at most 180 degrees 
with respect to the interior of $\hat{M}_i$.  Thus the conditions 
of Theorem 2.1 are satisfied along the line segment from 
$p_1$ to $p_2$.  
\end{proof}

\begin{proof}
({\sf of Theorem 1.5, the 2$n$-oids with alternating angles 
between the ends})

The proof of Theorem 1.5 is identical to the proof of Theorem 
1.1, once we replace the points $p_1$, ..., $p_7$ in the proof 
of Theorem 1.1 by the points 
$p_1 = (0,0,0)$, 
$p_2 = (0,-t\cos(\pi/n),-t\sin(\pi/n))$, 
$p_3 = (i,-t\cos(\pi/n),-t\sin(\pi/n))$, 
$p_4 = (i,-i,(t\cos(\pi/n) \newline -i)\tan(\theta/2)-t\sin(\pi/n))$, 
$p_5 = (-i,-i,(s-i)\tan(\theta/2))$, 
$p_6 = (-i,-s,0)$, and $p_7 = (0,-s,0)$.  
\end{proof}

\begin{proof}
({\sf of Theorem 1.6, classification})

Let $M$ be any complete genus-zero immersed catenoid-ended minimal 
surface with symmetry group $D_n \times \bfZ_2$ and at most 
$2n+1$ ends.  We can place $M$ in $\bfR^3$ so that its planes 
of reflectional symmetry are $P_0 = \{x_3 = 0\}$, 
$P_i = \{x_1 = \cot{\frac{i\pi}{n}}x_2\}$ for $i = 1,...,n-1$, and 
$P_n = \{x_2 = 0\}$.  Thus the $x_3$-axis is the axis for the 
rotational symmetry of order $n$ of the surface $M$.

Choose an orientation on $M$.  Consider an end 
$E$ of $M$ with limiting normal vector $\vec{v}$.  Let $l(E)$ be 
the central axis line of the end $E$.  Let Orb$(E)$ be the orbit of
$E$ under the symmetry group $D_n \times \bfZ_2$ of $M$.  

If $E$ is not invariant under any element of the symmetry 
group, then Orb$(E)$ would 
consist of $4n$ distinct ends, which contradicts our hypothesis.  
So $E$ must be invariant under reflection through $P_i$ for some $i$.
It follows that $\vec{v} \in 
P_i$ for some $i$.  Clearly $E$ cannot
be invariant under reflection through all $P_i$, so Orb$(E)$ 
must contain at least two ends.  In fact, the following list
represents all possibilities for Orb$(E)$: 
\begin{list}%
{\arabic{num})}{\usecounter{num}\setlength{\rightmargin}{\leftmargin}}

\item If the limiting normal vector $\vec{v}$ of $E$ is neither 
vertical nor horizontal, then Orb$(E)$ consists of 2$n$ ends.

\item If $\vec{v} \in P_0$ but $\vec{v} \not\in P_i$ for all $i \geq 
1$, then Orb$(E)$ consists of 2$n$ ends.

\item If $\vec{v} \in P_0$ and $\vec{v} \in P_i$ for some $i \geq 
1$ but the central axis $l(E) \not\in P_i$ for that value of $i \geq
1$, then Orb$(E)$ consists of 2$n$ ends.

\item If $\vec{v} \in P_0$ and $\vec{v} \in P_i$ for some $i \geq 
1$ and $l(E) \in P_i$ for that value of $i \geq 
1$, then Orb$(E)$ consists of $n$ ends.

\item If $\vec{v}$ is vertical and $l(E)$ is not the $x_3$-axis, 
then Orb$(E)$ consists of 2$n$ ends.

\item If $\vec{v}$ is vertical and $l(E)$ is the $x_3$-axis, 
then Orb$(E)$ consists of 2 ends.
\end{list}

Claim 1: $E$ is the unique end of $M$ with normal 
vector $v$ and axis $l(E)$.  

Suppose that $E^\prime$ is another end of $M$ with 
normal vector $v$, and that $l(E) = l(E^\prime)$.  
In all six cases above, the only end contained in Orb$(E)$ that has 
both the same normal vector and same central axis as $E$ is
$E$ itself.  Therefore, $E^\prime \not\in$ Orb$(E)$.  

If Orb$(E)$ contains $2n$ ends, then Orb$(E) \cup$Orb$(E^\prime)$ 
contains at least $2n+2$ ends, a contradiction.  
So in cases 1, 2, 3, and 5 above, $E^\prime$ cannot exist.  

Consider case 6.  In this case Orb$(E) \cup E^\prime$ consists of
three ends all with central axis the $x_3$-axis.  Thus all three of
these ends are invariant by a nontrivial rotation about the
$x_3$-axis.  This is impossible by Lemma~2.4, so $E^\prime$ cannot
exist.

Consider case 4.  In this case we may assume (by rotating in $\bfR^3$ 
about the $x_3$-axis if necessary) that $l(E)$ is the $x_1$-axis and
$\vec{v} = (1,0,0)$.  Let $R$ be rotation by 180 degrees about the
$x_1$-axis.  Note that Orb$(E) \cup$Orb$(E^\prime)$ contains all the
ends of $M$.  By Osserman's inequality, since we have a genus zero
minimal surface with 2$n$ embedded ends, we see that the Gauss map 
is a branched covering from $M$ to the unit sphere with 
order 2$n$-1.  Let $S$ be the set of 
2$n$-1 points on $M$ with Gauss map $(1,0,0)$ (including ends, and
counting with multiplicity).  $S$ is invariant under 
$R$ as a set, and there are an odd number of points in $S$.  
Since $R \circ R$ is the identity map, it follows that there must be 
an odd number of fixed points of $R$ contained in $S$.  The ends 
$E$ and $E^\prime$ represent two such points, so there must be a 
third.  But $R$ cannot have three fixed points, by Lemma~2.4.  
Thus $E^\prime$ cannot exist.  

This proves Claim 1.  

Claim 2: The ends of $M$ can have at most one orbit 
consisting of exactly two ends.  

Suppose $\{E_1,E_2\}$ is one orbit consisting of two ends, and suppose
$\{E_3,E_4\}$ is another orbit consisting of two ends.  Then 
$E_1$, $E_2$, $E_3$, and $E_4$ all have central axis the $x_3$-axis.  
Thus two of these ends must have the same normal vector, 
but that contradicts Claim 1.  

This proves Claim 2.  

So the only possibilities are the following:  
\begin{list}%
{\arabic{num})}{\usecounter{num}\setlength{\rightmargin}{\leftmargin}}

\item The ends of $M$ have a single orbit consisting of 
$n$ ends.  Then $M$ is equivalent to ${\cal J}{\cal M}_0(n)$.  

\item The ends of $M$ have two orbits, each consisting of 
$n$ ends.  Then $M$ is equivalent to 
${\cal A}{\cal W}_0(2n,w)$ with $w \neq 0,1$.    

\item The ends of $M$ have two orbits, one consisting of 
two ends, and one consisting of $n$ ends.  
Then $M$ is equivalent to ${\cal J}{\cal M}{\cal V}_0(n+2,w)$.  

\item The ends of $M$ have a single orbit consisting of 
$2n$ ends.  Then $M$ is equivalent to ${\cal P}_0(2n,\theta)$ 
or ${\cal A}{\cal A}_0(2n,\theta)$.  
\end{list}

$M$ cannot have only two ends, for in this case $M$ would be a 
catenoid (cf. \cite{Scn1}), whose symmetry group is not 
$D_n \times \bfZ_2$.  This completes the proof.  
\end{proof}

\section{Open problems}

One can ask whether the following immersed minimal surfaces with 
catenoid ends and symmetry group $D_n \times \bfZ_2$ exist:
\begin{list}%
{\arabic{num})}{\usecounter{num}\setlength{\rightmargin}{\leftmargin}}
\item The prismoids ${\cal P}_0(2n,\theta)$ and 
${\cal P}_{n-1}(2n,\theta)$ plus two 
vertical ends (analogous to the way 
${\cal J}{\cal M}{\cal V}_0(n+2,w)$ 
is the Jorge-Meeks 
surface plus two vertical ends).  
\item The surface ${\cal A}{\cal W}_0(2n,w)$ plus two vertical ends.  
\item A $3n$-oid with weight one at every third end as one 
travels around the circle of ends, and weight $w$ at 
the other $2n$ ends, $w \in (0,\infty)$.  
\item The example mentioned in the previous item, plus two vertical 
ends (assuming the surface is placed so that the first $3n$ ends 
have horizontal normal vectors).  
\item Prismoids with $k$ layers of catenoid ends, still with symmetry 
group $D_n \times \bfZ_2$, both genus zero and higher genus.  
If $k$ is odd, this 
surface would have $n$ ends with horizontal normal vectors, otherwise 
it would have no ends with horizontal normal vectors.  In either case,
it would have $n$ 
ends with normal vectors pointing upward making an angle $\theta_1$
with a horizontal plane, and $n$
ends with normal vectors pointing downward making an angle $\theta_1$
with a horizontal plane, $0 < \theta_1 < \pi/2$.  
The same would then be true for some angle $\theta_2$ with $\theta_1 
< \theta_2 < \pi/2$.  And again this holds for some angle 
$\theta_3$ with $\theta_2 
< \theta_3 < \pi/2$.  This continues up to the angle 
$\theta_{[k/2]}$, where $[k/2]$ is the greatest 
integer less than or equal to $k/2$.  
\item Prismoids with $k$ layers of catenoid ends plus two 
vertical ends.  
\end{list}

Solving some of the conjectures above might lead to a generalization 
of Theorem~1.6 to higher numbers of ends.  

A broader open question is the following:  

\begin{conjecture}
{\sf (Kusner's conjecture)}  
Any balanced configuration $\{v_1, ... ,v_n\}$ of $n$ vectors, 
such that for all $i$ and $j$, $v_i \neq r \cdot v_j$ for any
positive real $r$, can be realized as a 
genus-zero immersed minimal surface with $n$ catenoid ends.  
\end{conjecture}

Kapouleas \cite{Kap} has some corresponding results for the 
nonminimal constant-mean-curvature case.  

Shin Kato, Masaaki Umehara, and Kotaro Yamada have some results 
in the direction of these open questions 
\cite{KUY}, \cite{Kat}, \cite{UmYa}.


\begin{thebibliography}{KKMS}

\bibitem[An]{An} {\sc M. Anderson}, Curvature estimates for minimal 
surfaces in 3-manifolds, Ann. Sci. \'{E}cole Norm. Sup. 18 
(1985), 89-105. 
\bibitem[BeRo]{BeRo} {\sc J. Berglund and W. Rossman}, Minimal surfaces with 
catenoid ends, to appear in Pacific J. Math.  
\bibitem[HoMe]{HoMe} {\sc D. Hoffman and W. H. Meeks III}, Minimal surfaces 
based on the catenoid, 
Amer. Math. Monthly, Special Geometry Issue 97(8) (1990), 702-730.
\bibitem[JeSe]{JeSe} {\sc H. Jenkins and J. Serrin}, Variational problems of
minimal surface type II. Boundary value 
problems for the minimal surface equation, 
Arch. Rational Mech. Anal. 21 (1966), 321-342.  
\bibitem[JoMe]{JoMe} {\sc L. P. M. Jorge and W. H. Meeks III}, The topology
of complete minimal surfaces of finite total Gaussian curvature, 
Topology 22(2) (1983), 203-221.  
\bibitem[Kap]{Kap} {\sc N. Kapouleas}, Complete constant mean curvature 
surfaces in Euclidean three space, Ann. of Math. 131 (1990), 239-330.  
\bibitem[Ka1]{Ka1} {\sc H. Karcher}, Embedded minimal surfaces derived from
Scherk's surfaces, Manuscripta Math. 62 (1983), 83-114.
\bibitem[Ka2]{Ka2} {\sc H. Karcher}, The triply periodic minimal surfaces
of Alan Schoen and their constant mean curvature companions, 
Manuscripta Math. 64 (1989), 291-357.
\bibitem[Ka3]{Ka3} {\sc H. Karcher}, Construction of minimal surfaces,
Surveys in Geometry, 1-96, Univ. of Tokyo (1989) (Also:
Lecture Notes No. 12, SFB256, Bonn, 1989).
\bibitem[Ka4]{Ka4} {\sc H. Karcher}, Construction of higher genus embedded
minimal surfaces, preprint.
\bibitem[Kat]{Kat} {\sc S. Kato}, Construction of $n$-end catenoids 
with prescribed flux, preprint.
\bibitem[KKS]{KKS} {\sc N. Korevaar, R. Kusner, and B. Solomon}, The 
structure of complete embedded surfaces with constant mean curvature, 
J. Differential Geom. 30 (1989), 465-503.
\bibitem[KUY]{KUY} {\sc S. Kato, M. Umehara, and K. Yamada}, personal 
communication.  
\bibitem[MESH]{MESH} {\sc J. Hoffman}, Software for constructing minimal
surfaces using Weierstrass data, G.A.N.G. Lab, Univ. of
Massachusetts, Amherst (1980-1993).
\bibitem[MeYa]{MeYa} {\sc W. Meeks and S-T Yau}, The existence of embedded 
minimal surfaces and the problem of uniqueness, 
Math. Z. 179 (1982), 151-168.
\bibitem[Ni]{Ni} {\sc J. C. C. Nitsche}, 
Uber ein verallgemeinertes Dirichletsches
Problem fur die Minimalflachengleichung und hebbare Unstetigkeiten
ihrer
Losungen, Math. Ann. 158 (1965), 302-214.
\bibitem[Os]{Os} {\sc R. Osserman}, A survey of minimal surfaces, 
Dover, New York, 1986.
\bibitem[Scn1]{Scn1} {\sc R. Schoen}, Uniqueness, symmetry, and embeddedness
of minimal surfaces, J. Differential Geom. 18 (1982), 791-809.
\bibitem[Scn2]{Scn2} {\sc R. Schoen}, Estimates for stable minimal surfaces in 
three dimensional manifolds, 
Ann. of Math. Studies, 103, Princeton Univ. Press, 1983.
\bibitem[UmYa]{UmYa} {\sc M. Umehara and K. Yamada}, Surfaces of constant 
mean curvature $c$ in $\bfH^3(-c^2)$ with prescribed hyperbolic 
Gauss map, preprint.  
\bibitem[Xu]{Xu} {\sc Y. Xu}, Symmetric minimal surfaces in $\bfR^3$, 
to appear in Pacific J. Math.  

\vspace{1.31in}

{\sc Mathematical Institute}

{\sc Faculty of Science}

{\sc Tohoku University}

{\sc Sendai 980}

{\sc Japan}

\end{thebibliography}
\end{document}